\documentclass[aos, preprint]{imsart}
\setattribute{journal}{name}{}

\usepackage{geometry}
\usepackage{amsfonts}
\usepackage{amsmath}
\usepackage{amssymb}
\usepackage{amsthm}
\usepackage{mathrsfs}
\usepackage[numbers]{natbib}
\usepackage[usenames, dvipsnames]{color}
\usepackage{hyperref}
\usepackage{graphicx}
\usepackage{caption}
\usepackage{float}

\usepackage[T1]{fontenc}
\usepackage[utf8]{inputenc}

\startlocaldefs

\renewcommand\labelenumi{(\roman{enumi})}
\renewcommand\theenumi\labelenumi
\newcommand{\Ex}{{\mathbb{E}}}

\newcommand{\Cov}{\mathbb{C}\text{ov}}

\newcommand{\summ}{{\sum_{m=-\ell}^{\ell}}}

\newcommand{\op}{\text{op}}
\newcommand{\tr}{\text{Tr}}
\newtheorem{theorem}{Theorem}

\newtheorem{condition}[theorem]{Condition}

\newtheorem{corollary}[theorem]{Corollary}

\newtheorem{definition}[theorem]{Definition}
\newtheorem{example}[theorem]{Example}

\newtheorem{lemma}[theorem]{Lemma}

\newtheorem{remark}[theorem]{Remark}

\endlocaldefs

\input{tcilatex}

\begin{document}

\begin{frontmatter}

\title{Asymptotics for Spherical Functional Autoregressions}
\runtitle{Spherical Functional Autoregressions}


\author{\fnms{Alessia} \snm{Caponera}\ead[label=e1]{alessia.caponera@uniroma1.it}\thanksref{t1} and
\fnms{Domenico} \snm{Marinucci}\ead[label=e2]{marinucc@mat.uniroma2.it}\thanksref{t1}}

\address{Piazzale Aldo Moro, 5\\ 00185 Rome \\Italy\\\href{mailto:alessia.caponera@uniroma1.it}{E-mail: \textnormal{alessia.caponera@uniroma1.it}}
}
\affiliation{Department of Statistical Sciences, Sapienza University of Rome}


\thankstext{t1}{DM acknowledges the MIUR Excellence Department Project awarded to the Department of Mathematics, University of Rome Tor Vergata, CUP E83C18000100006. We are also grateful to Pierpaolo Brutti for many insightful suggestions and conversations.}
\address{Via della Ricerca Scientifica, 1\\00133 Rome \\Italy\\ \href{mailto:marinucc@mat.uniroma2.it}{E-mail: \textnormal{marinucc@mat.uniroma2.it}}
}
\affiliation{Department of Mathematics, University of Rome Tor Vergata}

\runauthor{Caponera and Marinucci}

\begin{abstract}
In this paper, we investigate a class of spherical functional autoregressive
processes, and we discuss the estimation of the corresponding autoregressive
kernels. In particular, we first establish a consistency result (in sup and
mean-square norm), then a quantitative central limit theorem (in Wasserstein
distance), and finally a weak convergence result, under more restrictive
regularity conditions. Our results are validated by a small numerical
investigation.
\end{abstract}

\begin{keyword}[class=MSC]
\kwd[Primary ]{62M15}
\kwd[; secondary ]{60G15 60F05 62M40}
\end{keyword}

\begin{keyword}
\kwd{Spherical Functional Autoregressions}
\kwd{Spherical Harmonics}
\kwd{Quantitative Central Limit Theorem}
\kwd{Wasserstein Distance}
\kwd{Weak Convergence}
\end{keyword}

\end{frontmatter}


\section{Introduction}

In recent years, a lot of interest has been drawn by the statistical
analysis of spherical isotropic random fields. These investigations have
been motivated by a wide array of applications arising in many different
areas, including in particular, Cosmology, Astrophysics, Geophysics, Climate
and Atmospheric Sciences, and many others, see, e.g., \cite{BKMP,Berg,Clarke,CM2018,Gneiting,Fan,Fan2,Leonenko,MarVad,Porcu}. Most
papers in Cosmology and Astrophysics have focussed so far on spherical
random fields with no temporal dependence; the next generation of
Cosmological experiments
is however going to make the time dependence much more relevant. On the other hand, applications in Climate,
Atmospheric Sciences, Geophysics, and several other areas have always been
naturally modelled in terms of a double-dependence in the spatial and
temporal domains. In many works of these fields, the attention has been
focussed on the definition of wide classes of space-time covariance functions,
and then on the derivation of likelihood functions; the literature on these
themes is vast and we make no attempt to a complete list of references, see
for instance \cite{Clarke, Gneiting, Jun, Porcu} and the references therein.

Our purpose in this paper is to investigate a class of space-time processes,
which can be viewed as functional autoregressions taking values in $L^{2}({\mathbb{S}^{2}})$; we refer to \cite{Bosq} for a general textbook analysis
of functional autoregressions taking values in Hilbert spaces, and \cite%
{Aue,Hormann,Panaretos} for a very partial list of some important recent
references.

Dealing with functional spherical autoregressions ensures some very convenient simplifications; in
particular, we exploit the analytic properties of the standard orthonormal basis of $L^{2}({\mathbb{S}^{2}})$ and some natural isotropy requirements to obtain
neat expressions for the autoregressive functionals, which are then
estimated by a form of frequency-domain least squares. For our estimators,
we are able to establish rates of consistency (in $L^{2}$ and $L^{\infty }$
norms) and a quantitative version of the Central Limit Theorem, in
Wasserstein distance. In particular, we derive explicit bounds for the rate
of convergence to the limiting Gaussian distribution by means of the rich
machinery of Stein-Malliavin methods (see \cite{noupebook}); to the
best of our knowledge, this is the first Quantitative Central Limit Theorem
established in the framework of functional-valued stationary processes.
Under stronger regularity conditions, we are able to establish a weak
convergence result for the kernel estimators; our results are then
illustrated by simulations.

The plan of our work is then as follows: in Section \ref{BackgroundNotation}
we present background results on the harmonic analysis of spherical random
fields and on Stein-Malliavin methods. In Section \ref{SpaceTime} we present
our basic model; we show how, under isotropy, the model enjoys a number of symmetry
properties which greatly simplify our approach. Our main results are then
collected in Section \ref{MainResults}, where we investigate rates of
convergence and the Quantitative Central Limit Theorem; we consider also
weak convergence in $C_{p}\left( [-1,1]\right) $, under stronger regularity
conditions for the autoregressive kernels. Large parts of the proofs and
many auxiliary lemmas, some of possible independent interest, are collected
in Sections \ref{Proofs} and in the Appendix (Supplementary Material). 
Finally, Section \ref{Numerical} provides numerical estimates on the behaviour of our procedures.

\section{Background and Notation} \label{BackgroundNotation}

\subsection{Spectral Representation of Isotropic Random Fields on the Sphere}

Let $\{T(x),\ x\in {\mathbb{S}^{2}}\}$ denote a finite variance, isotropic
random field on the unit sphere ${\mathbb{S}^{2}}=\{x\in \mathbb{R}%
^{3}:\Vert x\Vert =1\}$, by which we mean as usual that $T(g\cdot )\overset{d%
}{=}T(\cdot ),\ \forall g\in SO(3)$ the standard 3-dimensional group of
rotations; here the identity in distribution must be understood in
the sense of stochastic processes: for notational simplicity, and without
loss of generality, we will assume in the sequel that ${\mathbb{E}}[T(x)]=0$. It is well-known that the following representation holds, in the
mean-square sense:
\begin{equation}
T(x)=\sum_{\ell =0}^{\infty }T_{\ell }(x)\ ,\qquad T_{\ell }(x)=\sum_{m=-\ell
}^{\ell }a_{\ell ,m}Y_{\ell ,m}(x) \text{ ,}  \label{spectral_representation}
\end{equation}%
where $\{Y_{\ell ,m}(\cdot ),\ \ell \geq 0,\ -\ell \leq m\leq \ell \}$ are
the standard basis of spherical harmonics, which satisfy (for $\varphi \in
\lbrack 0,2\pi )$\ , $\vartheta \in 0,\pi ]$)
\begin{equation*}
\Delta _{{\mathbb{S}^{2}}}Y_{\ell ,m}=-\ell (\ell +1)Y_{\ell ,m}\text{ ,
} \quad \Delta _{\mathbb{S}^{2}}:=\frac{1}{\sin \vartheta }\frac{\partial }{\partial \vartheta }\left( \sin \vartheta \frac{\partial }{\partial
\vartheta }\right) +\frac{1}{\sin ^{2}\vartheta }\frac{\partial }{\partial
\varphi ^{2}}\text{ ; }
\end{equation*}%
also $\{a_{\ell ,m},\ \ell \geq 0,\ -\ell \leq m\leq \ell \}$ are a
triangular array of zero-mean, real-valued random coefficients whose
covariance structure is given by
\begin{equation*}
{\mathbb{E}}[a_{\ell ,m}a_{\ell ^{\prime },m^{\prime }}]=C_{\ell }\delta
_{\ell }^{\ell ^{\prime }}\delta _{m}^{m^{\prime }};
\end{equation*}%
here $\delta _{a}^{b}$ is the Kronecker delta function, and the sequence $%
\{C_{\ell },\ \ell \geq 0\}$ represents the angular power spectrum of the
field. Throughout this paper we consider the real-valued basis of spherical
harmonics, and therefore the random coefficients are real-valued random
variables for all $(\ell ,m)$ (we refer for instance to \cite{MaPeCUP} for a
more detailed discussion on spectral representations on the sphere). Note
that the random coefficients $\left\{ a_{\ell ,m}\right\} $ can be obtained
by a direct inversion formula from the map $T(\cdot ),$ indeed we have%
\begin{equation*}
a_{\ell ,m}:=\int_{{\mathbb{S}^{2}}}T(x)Y_{\ell ,m}(x)dx\text{ }.
\end{equation*}%
Here, we recall also the following \emph{addition formula} for spherical
harmonics (see \cite{MaPeCUP}, equation 3.42) which entails that, for any $%
x,y\in {\mathbb{S}^{2}}$,
\begin{equation}
{\sum_{m=-\ell }^{\ell }}Y_{\ell ,m}(x)Y_{\ell ,m}(y)=\frac{2\ell +1}{4\pi }%
P_{\ell }(\langle x,y\rangle )\ ,  \label{addition}
\end{equation}%
where $\langle x,y\rangle $ denotes the standard inner product in $\mathbb{R}^3$, and $P_{\ell
}(\cdot )$ represents the $\ell $-th Legendre polynomial, defined as usual
by
\begin{equation*}
P_{\ell }(t)=\frac{1}{2^{\ell }\ell !}\frac{d^{\ell }}{dt^{\ell }}%
(t^{2}-1)^{\ell },\qquad t\in \lbrack -1,1]\, , \ \ell \geq 0 \ .
\end{equation*}%
It is easy to show that $P_{\ell }(1)=1;$ moreover, the following \textit{%
reproducing property} is satisfied, i.e.,
\begin{equation*}
\int_{{\mathbb{S}^{2}}}\frac{2\ell +1}{4\pi }P_{\ell }(\langle x,y\rangle )%
\frac{2\ell +1}{4\pi }P_{\ell }(\langle y,z\rangle )dy=\frac{2\ell +1}{4\pi }%
P_{\ell }(\langle x,z\rangle )\ .
\end{equation*}%
Under isotropy, from \eqref{spectral_representation} and \eqref{addition} the
covariance function $\Gamma (x,y)={\mathbb{E}}[T(x)T(y)]$ satisfies
\begin{align*}
\Gamma (x,y)&={\sum_{\ell =0}^{\infty }\sum_{m=-\ell }^{\ell }}C_{\ell
}Y_{\ell ,m}(x)Y_{\ell ,m}(y)\\
&=\sum_{\ell =0}^{\infty }C_{\ell }\frac{2\ell +1%
}{4\pi }P_{\ell }(\langle x,y\rangle ) 	\ , \qquad \text{for all } x, y\in \mathbb{S}^{2} \ .
\end{align*}

In the sequel, given any two positive sequences $\{a_k,\  k\in \mathbb{N}\}$, $\{b_k, \ k\in \mathbb{N}\}$, we shall write $a_k\sim b_k$ if $\exists
c_{1},c_{2}>0$ such that $c_{1}b_k\leq a_k\leq c_{2}b_k,\forall k\in \mathbb{N}$. In addition, we will denote with $const$ a positive real constant,
which may change from line to line; also, we use $\|\cdot\|_{L^2(\mathbb{S}^2)}$ for the usual $L^2$ norm on the sphere, $\Lambda_{\min}(A)$ and $\Lambda_{\max}(A)$ for the minimum and maximum eigenvalues of the matrix $A$, respectively, $\Vert A\Vert _{%
\text{op}}$ for the operator norm of $A$, i.e., $\Vert A\Vert _{\text{op%
}}\ =\sqrt{\lambda _{\max }(A^{\prime }A)}$,
and $\text{Tr}(A)$ for the
trace of $A$.

\subsection{Hermite Polynomials and Stein-Malliavin Results}

Let us recall the family of Hermite polynomials $%
\{H_{q}(\cdot ),\ q\geq 0\}$, defined by
\begin{equation*}
H_{q}(x):=(-1)^{q}e^{x^{2}/2}\frac{d^{q}}{dx^{q}}e^{-x^{2}/2}\ ,\qquad x\in
\mathbb{R}\text{ };
\end{equation*}%
for instance, the first few are given by $H_{1}(x)=x$, $H_{2}(x)=x^{2}-1$, $%
H_{3}(x)=x^{3}-3x$ and $H_{4}(x)=x^{4}-6x^{2}+3$. The sequence $%
\{(q!)^{-1/2}H_{q}(\cdot ),\ q\geq 0\}$ is an orthonormal basis of the space
of finite variance transform of Gaussian variables, i.e.,%
\begin{equation*}
G(X)=\sum_{q=0}^{\infty }J_{q}(G)\frac{H_{q}(X)}{q!}\ , \quad \text{for }G\text{
s.t. }\mathbb{E}G^{2}(X)<\infty\ , \quad J_{q}(G):=\mathbb{E}%
G(X)H_{q}(X)\text{ };
\end{equation*}%
more generally, for the space $L^{2}(\Omega)$ generated by any Gaussian random measure we can write the Stroock-Varadhan decomposition
\begin{equation*}
L^{2}(\Omega )=\bigoplus\limits_{q=0}^{\infty }\mathcal{H}_{q}\text{ ,}
\end{equation*}
where $\mathcal{H}_{q}$ is the $q$-th order Wiener chaos, i.e., the space spanned by linear combinations of $q$-th order Hermite polynomials,
see \cite{noupebook} for more discussions and details.

We shall exploit extensively a very powerful technique, recently discovered
by \cite{noupe2009}, to establish Quantitative Central Limit Theorems for
sequences of random variables belonging to Wiener chaoses. To explain what
we mean by a Quantitative Central Limit Theorem, we recall first the notion
of \emph{Wasserstein distance}, i.e., for any two $d$-dimensional random
variables $X,Y,$
\begin{equation*}
d_{W}(X,Y)=\sup_{h(\cdot):\Vert h\Vert _{\text{Lip}}\leq 1}|{\mathbb{E}}[h(X)]-{%
\mathbb{E}}[h(Y)]| \ , \quad \text{where }\Vert h\Vert _{\text{Lip}}=\sup_{\substack{ %
x\neq y  \\ x,y\in \mathbb{R}^{d}}}\frac{|h(x)-h(y)|}{\Vert x-y\Vert }\ ,
\end{equation*}%
with $\Vert \cdot \Vert $ the usual Euclidean norm on $\mathbb{R}^{d}$,
where we assume that ${\mathbb{E}}|h(X)|<\infty $, ${\mathbb{E}}%
|h(Y)|<\infty $ for every $h(\cdot)$. See \cite{noupebook} for a discussion of the
main properties of $d_{W}$ and for other examples of probability metrics;
here we recall simply that
\begin{equation}
d_{W}(X,Y)\leq {\mathbb{E}}\Vert X-Y\Vert \text{ .}  \label{E_dist}
\end{equation}

It is shown in \cite{noupebook} that for sequences of zero-mean scalar
random variables $\left\{ Z_k,\ k \in \mathbb{N} \right\}$ belonging to $%
\mathcal{H}_{q}$ ($q\geq 2$) such that ${\mathbb{E}}[Z_k]=\sigma ^{2}>0$,
one has the remarkable inequality%
\begin{equation}
d_{W}(Z,Z_k)\leq \frac{1}{\sigma }\sqrt{\frac{2q-2}{3\pi q}\big(\mathbb{E}%
[Z_k^{4}]-3\sigma ^{4}\big)}\ ,  \label{npbound}
\end{equation}%
where $Z\overset{d}{=}\mathcal{N}(0,\sigma ^{2})$ (in our proof below we
will actually exploit a multivariate extension of this inequality, also
given in \cite{noupebook}). The inequality in \eqref{npbound} can be proved by
means of the so-called Stein-Malliavin approach, which establishes a deep
and surprising connection between Malliavin calculus and Stein's equation as
a tool for the investigation of limiting distributions. In particular, in
view of \eqref{npbound} for sequences that belong to Wiener chaoses the
investigation of the asymptotic behaviour of the fourth-moment is enough to
investigate not only the validity of a central limit theorem, but also the
rate of convergence to the Gaussian limiting distribution.

\section{Spherical Random Fields with Temporal Dependence}  \label{SpaceTime}

We are now ready to introduce our model of interest. As usual, by space-time
spherical random fields we mean a collection of random variables $\{T(x,t),\
(x,t)\in {\mathbb{S}^{2}}\times \mathbb{Z}\}$ such that the application $%
T:\Omega \times {\mathbb{S}^{2}}\times \mathbb{Z}\rightarrow \mathbb{R}$ is $%
\Im\otimes \mathcal{B}({\mathbb{S}^{2}}\times \mathbb{Z})$%
-measurable, for some probability space $(\Omega ,\Im,\mathbb{P})$.
The following definition is standard:

\begin{definition}
$\{T(x,t),\ (x,t)\in {\mathbb{S}^{2}}\times \mathbb{Z}\}$ is 2-weakly
isotropic stationary if ${\mathbb{E}}[T(x,t)]$ is constant $\forall (x,t)\in
{\mathbb{S}^{2}}\times \mathbb{Z}$ and the covariance function $\Gamma $ is
a spatially isotropic and temporally stationary function on $({\mathbb{S}^{2}%
}\times \mathbb{Z})^{2}$, that is there exists $\Gamma _{0}:[-1,1]\times
\mathbb{Z}\rightarrow \mathbb{R}$ such that
\begin{equation*}
\Gamma (x,t,y,s)=\Gamma _{0}(\langle x,y\rangle ,t-s)\text{ },\qquad
\forall (x,t),(y,s)\in {\mathbb{S}^{2}}\times \mathbb{Z}\ .
\end{equation*}
\end{definition}

In particular, we will focus on Gaussian random fields, where of course weak
stationarity entails strong isotropy and stationarity, i.e., the law of $%
T(g\cdot ,\cdot +\tau )$ is the same as the law of $T(\cdot ,\cdot ),$ in
the sense of processes, for all $g\in SO(3)$ and $\tau \in \mathbb{Z}$. Note
that, for (zero-mean) finite variance random fields, $T_{t}(\cdot )\equiv
T(\cdot ,t)$ is a random function of $L^{2}({\mathbb{S}^{2}})$, $t\in
\mathbb{Z}$. Thus, for any fixed $t\in \mathbb{Z}$, the following spectral
representation holds:
\begin{equation}
T_{t}(x)=\sum_{\ell =0}^{\infty }\sum_{m=-\ell }^{\ell }a_{\ell
,m}(t)Y_{\ell ,m}(x)\ ,  \label{spectral_representation_space-time}
\end{equation}%
where $\{Y_{\ell ,m}(\cdot ),\ \ell \geq 0,\ -\ell \leq m\leq \ell \}$ are
spherical harmonics, and $\{a_{\ell ,m}(t),\ \ell \geq 0,\ -\ell \leq m\leq
\ell \}$ (zero-mean) random coefficients which satisfy
\begin{equation*}
\qquad {\mathbb{E}}[a_{\ell ,m}(t)a_{\ell ^{\prime },m^{\prime
}}(s)]=C_{\ell }(t-s)\delta _{\ell }^{\ell ^{\prime }}\delta _{m}^{m^{\prime
}}\ ,\qquad t,s\in \mathbb{Z}\ .
\end{equation*}%
Note that $\{C_{\ell }(0),\,\ell \geq 0\}$ corresponds to the angular power
spectrum of the spherical field at a given time point, for which we will
simply write $\{C_{\ell }\}$. As for the isotropic case, for fixed $t,s\in
\mathbb{Z}$, the covariance function $\Gamma (x,t,y,s)$ is easily shown to
have a spectral decomposition in terms of Legendre polynomials (Schoenberg's
Theorem, see also \cite{Berg}), i.e., for every $(x,t),(y,s)\in {\mathbb{S}^{2}}\times \mathbb{Z}$,
\begin{equation*}
\Gamma (x,t,y,s)=\sum_{\ell =0}^{\infty }C_{\ell }(t-s)\frac{2\ell +1}{4\pi }%
P_{\ell }(\langle x,y\rangle )\text{ .}
\end{equation*}

\begin{remark}
By exploiting results from \cite{Panaretos}, it would also be possible to
rewrite \eqref{spectral_representation_space-time} by means of the
Cram\'er-Karhunen-Lo\'eve representation%
\begin{equation*}
T_{t}(x)=\int_{-\pi }^{\pi }\sum_{\ell =0}^{\infty }\sum_{m=-\ell }^{\ell
}\exp (-i\lambda t)Y_{\ell ,m}(x)dW_{\ell ,m}(\lambda )\ , \quad \text{in }%
L^{2}(\Omega )\text{ ,}
\end{equation*}%
where $\left\{W_{\ell ,m}(\cdot)\right\} $ is a family of independent
complex-valued Gaussian random measures, with mean zero and covariance
structure%
\begin{equation*}
\mathbb{E}\left[ W_{\ell ,m}(A)\overline{W}_{\ell ,m}(B)\right] =\int_{A\cap
B}f_{\ell }(\lambda )d\lambda \ , \quad \text{for all }A,B\subset \lbrack -\pi ,\pi
]\text{ ,}
\end{equation*}%
where $ f_{\ell }(\cdot)$ denotes the spectral density of the
process $\left\{ a_{\ell ,m}(t), \ t \in \mathbb{Z} \right\} ,$ which is introduced below and
satisfies%
\begin{equation*}
\mathbb{E}\left[ a_{\ell ,m}(t)a_{\ell ,m}(t+\tau )\right] =\int_{-\pi
}^{\pi }\exp (i\lambda \tau )f_{\ell }(\lambda )d\lambda \text{ .}
\end{equation*}%
This approach is not pursued here, see also \cite{Caponera} for more discussion
and details.
\end{remark}

\subsection{Spherical Autoregressions}

Consider now two zero-mean
Gaussian isotropic stationary random fields
\[
\{T(x,t),\ (x,t)\in {\mathbb{S}^{2}}\times \mathbb{Z}\} \text{ and }\{Z(x,t),\ (x,t)\in {\mathbb{S}^{2}}\times \mathbb{Z}\} \ ,
\]
so that ${\mathbb{E}}[T^{2}(x,t)]<\infty $ and ${\mathbb{E}}[Z^{2}(x,t)]<\infty $. Let us start from the definition of a Gaussian spherical white noise process.

\begin{definition}[Gaussian Spherical White Noise]
\label{SphericalWhiteNoise} $\{Z(x,t),\ (x,t)\in {\mathbb{S}^{2}}\times \mathbb{Z}\}$ is a sequence of independent and
identically distributed Gaussian isotropic spherical random fields. That is

a) for every fixed $t\in \mathbb{Z}$, $Z(\cdot,t)$ is a Gaussian, zero-mean
isotropic random field, with covariance function%
\begin{equation*}
\Gamma _{Z}(x,y)=\sum_{\ell =0}^{\infty }\frac{2\ell +1}{4\pi }C_{\ell
;Z}P_{\ell }(\left\langle x,y\right\rangle )\text{ ,} \quad \sum_{\ell
=0}^{\infty }\frac{2\ell +1}{4\pi }C_{\ell; Z}<\infty \text{ ; }
\end{equation*}%
here, $\left\{ C_{\ell; Z}\right\} $ denotes as usual the angular power
spectrum of $Z(\cdot,t)$;

b) for every $t\neq s,$ the random fields $Z(\cdot,t)$ are independent.
\end{definition}

\begin{remark}
Note that we are defining the field as a collection of random variables
defined on every pair $(x,t)\in S^{2}\times \mathbb{Z}.$ Alternatively, one
could view the fields as random elements in a Hilbert space (in our case,
corresponding to $L^{2}(\mathbb{S}^{2}),$ see \cite{Bosq}$);$ the two
approaches are equivalent here, because throughout this paper we will always
be dealing with mean-square continuous random fields.
\end{remark}

\begin{definition}
A spherical isotropic kernel operator is an application $\Phi
:L^{2}(S^{2})\rightarrow L^{2}(S^{2})$ which satisfies%
\begin{equation*}
(\Phi f)(\cdot)=\int_{S^{2}}k(\left\langle \cdot,y\right\rangle )f(y)dy\text{ ,}%
\qquad \text{ some }k(\cdot)\in L^{2}[-1,1]\text{ .}
\end{equation*}
\end{definition}

The following representation holds, in the $L^{2}$ sense, for the kernel
associated to $\Phi $:
\begin{equation}
k(\left\langle x,y\right\rangle )=\sum_{\ell =0}^{\infty }\phi _{\ell }\frac{%
2\ell +1}{4\pi }P_{\ell }(\left\langle x,y\right\rangle )\text{ .}
\label{kernel_expansion}
\end{equation}%
The coefficients $\{\phi _{\ell },\ \ell \geq 0\}$ corresponds to the eigenvalues of the
operator $\Phi $ and the associated eigenfunctions are the family of
spherical harmonics $\left\{ Y_{\ell ,m}\right\}$, yielding%
\begin{equation*}
\Phi Y_{\ell ,m}=\phi _{\ell }Y_{\ell ,m}\text{ ,}
\end{equation*}%
Thus, it holds $\sum_{\ell }(2\ell +1)\phi _{\ell }^{2}<\infty $, and hence
this operator is Hilbert-Schmidt (see, e.g., \cite{Hsing}). In this paper, we
shall also consider trace class operators, namely such that $\sum_{\ell
}(2\ell +1)|\phi _{\ell }|<\infty $, for which the representation %
\eqref{kernel_expansion} holds pointwise for every $x,y\in {\mathbb{S}^{2}}$.

\begin{definition}
\label{SPHAR(p)}The collection of random variables $\left\{ T(x,t), (x,t) \in \mathbb{S}^{2} \times \mathbb{Z}\right\}$ satisfies the Spherical
Autoregressive process of order $p$ (written $SPHAR(p)$) if there exist $p$
isotropic kernel operators $\left\{ \Phi _{1},\dots ,\Phi _{p}\right\} $ such that
\begin{equation}
T_{t}(x)-(\Phi _{1}T_{t-1})(x)-\cdots-(\Phi _{p}T_{t-p})(x)-Z_{t}(x)=0\text{
,}  \label{defar}
\end{equation}%
for all $(x,t)\in \mathbb{S}^{2}\times \mathbb{Z}$, the equality holding
both in the $L^{2}(\Omega )$ and in the $L^{2}(\Omega \times \mathbb{S}^{2})$
sense.
\end{definition}

\begin{remark}
It should be noted that the solution process is defined pointwise, i.e., for
each $(x,t)$ there exists a random variable defined on ($\Omega ,\Im ,%
\mathbb{P}$) such that the identity \eqref{defar} holds.
\end{remark}

Let us define the eigenvalues $\left\{ \phi _{\ell
;j}, \ \ell \ge 0, \ j=1,\dots,p\right\}$ which satisfy%
\begin{equation*}
\Phi _{j}Y_{\ell ,m}=\phi _{\ell ;j}Y_{\ell ,m}\text{ , } \text{ and hence }k_{j}(\langle x,y \rangle )=\sum_{\ell =0}^{\infty }\phi _{\ell ;j}\frac{%
2\ell +1}{4\pi }P_{\ell }(\langle x,y\rangle )\text{ .}
\end{equation*}%
Hence, for any $t\in \mathbb{Z}$,
\begin{equation*}
(\Phi _{j}T_{t-j})(x)=\sum_{\ell =0}^{\infty }\sum_{m=-\ell }^{\ell }\phi
_{\ell ;j}a_{\ell ,m}(t-j)Y_{\ell ,m}(x)\ ,
\end{equation*}%
that is $(\Phi _{j}T_{t-j})(\cdot )$ admits a spectral representation in
terms of spherical harmonics with coefficients $\{\phi _{\ell ;j}a_{\ell
m}(t-j),\ \ell \geq 0,\ -\ell \leq m\leq \ell \}$. \ Likewise, we obtain
\begin{equation}
a_{\ell ,m}(t)=\phi _{\ell ;1}a_{\ell ,m}(t-1)+\cdots +\phi _{\ell
;p}a_{\ell ,m}(t-p)+a_{\ell ,m;Z}(t)\text{ ;}  \label{arp}
\end{equation}%
to ensure identifiability, we assume that there exists at least an $\ell $
such that $\phi _{\ell ;p}\neq 0,$ so that $\Pr \{(\Phi _{p}T_{t})(\cdot
)\neq 0\}>0,\ t\in \mathbb{Z},$ see again \cite{Bosq}. Let us now define as
usual the associated polynomials $\phi _{\ell }:\mathbb{C\rightarrow C},\
\ell \geq 0$:
\begin{equation}
\phi _{\ell }(z)=1-\phi _{\ell ;1}z-\cdots -\phi _{\ell ;p}z^{p}\text{ .}
\label{AssoPoly}
\end{equation}

\begin{condition}
\label{stationarity} The sequence of polynomials \eqref{AssoPoly} is such
that $|z| \le 1 + \delta \ \Rightarrow \ \phi _{\ell }(z)\ne 0$, some $\delta>0$ . More explicitly,
there are no roots in a $\delta$-enlargement of the unit disk, for all $\ell \geq 0$.
\end{condition}

\begin{remark}
Under Condition \ref{stationarity}, Equation \eqref{defar} admits a unique stationary isotropic solution; the proof can be given along the same lines as in \cite{Bosq}, and it is omitted for brevity's sake, see \cite{Caponera} for more discussion and details.
\end{remark}

\begin{example}[$SPHAR(1)$] The family of random variables $\{T(x,t), \ (x, t) \in S^{2}\times \mathbb{Z} \}$ is a spherical $AR(1)$ process if for all pairs $(x,t)\in S^{2}\times Z$
it satisfies%
\begin{equation*}
T_{t}(x)=(\Phi _{1}T_{t-1})(x)+Z_{t}(x)\text{ .}
\end{equation*}%
In this case, the Condition \ref{stationarity} simply becomes $|\phi _{\ell
}|<\frac{1}{1+\delta}$, $\ell \geq 0$ .
\end{example}

\begin{remark}
The autocovariance function of a stationary spherical $AR(1)$ process is
easily seen to be given by (writing $\tau =t-s)$%
\begin{equation*}
\Gamma (x,t,y,s)=\Gamma _{0}(\langle x,y\rangle ,\tau )=\sum_{\ell
=0}^{\infty }C_{\ell }(\tau )\frac{2\ell +1}{4\pi }P_{\ell }(\langle
x,y\rangle )=\sum_{\ell =0}^{\infty }\frac{\phi _{\ell }^{|\tau |}C_{\ell ;Z}%
}{1-\phi _{\ell }^{2}}\frac{2\ell +1}{4\pi }P_{\ell }(\langle x,y\rangle )%
\text{ .}
\end{equation*}%
It is easy hence to envisage a number of parametric models for sphere-time
covariances; for instance, a simple proposal is%
\begin{eqnarray}
\phi _{\ell } &=&G\times \left\{ \left\vert \ell -\ell ^{\ast }\right\vert
+1\right\} ^{-\alpha _{\phi }}\text{ , }\ell ^{\ast }\ge 0 \text{ , }\alpha
_{\phi }>2\text{ , }0<G<1\text{ ,}  \label{Cupola} \\
C_{\ell ;Z} &=&G_{Z}(1+\ell )^{-\alpha _{Z}}\text{ , }\alpha _{Z}>2\text{ .}
\notag
\end{eqnarray}%
Here, the parameters $\alpha _{Z},\alpha _{\phi }$ control, respectively,
the smoothness of the innovation process and the regularity of the
autoregressive kernel (see \cite{LangSchwab}); the positive integer $\ell
^{\ast }$ can be seen as a sort of "characteristic scale", where the power
of the kernel is concentrated. More generally, we can take $\phi _{\ell
}=G(\ell ;\alpha _{1},\dots,\alpha _{q}),$ where $\alpha _{1},\dots,\alpha _{q}$
are fixed parameters and $G$ is any function such that%
\begin{equation*}
\sup_{\ell }\left\vert G(\ell ;\alpha _{1},\dots,\alpha _{q})\right\vert <1%
\text{ and }\sum_{\ell }(2\ell +1)\left\vert G(\ell ;\alpha _{1},\dots,\alpha
_{q})\right\vert <\infty \text{ , }
\end{equation*}%
uniformly over all values of $(\alpha _{1},\dots,\alpha _{q}).$
\end{remark}

\begin{condition}[Identifiability]
\label{Identifiability} The Gaussian spherical white noise
process $\left\{ Z(x,t)\right\} $ is such that $C_{\ell ;Z}>0$ for all $\ell
=0,1,2,\dots.$
\end{condition}

\begin{remark}
The previous condition is an identifiability assumption; indeed, it is
simple to verify from our arguments below that for $C_{\ell ;Z}=0$ the
component of the kernel corresponding to the $\ell $-th multipole is not
observable, i.e., the $AR(p)$ process has the same distribution whatever the
value of $\phi _{\ell }.$ It is possible, however, to estimate the
"sufficient" version of the kernel, i.e., its projection on the relevant
subspace, such that $C_{\ell ,Z}>0.$ The extension is straightforward and we
avoid it just for brevity and notational simplicity. Of course, as a
consequence we have that
\begin{equation*}
\int_{\mathbb{S}^{2}}\Gamma _{Z}(x,y)f(x)f(y)dxdy>0\ ,\qquad \forall f(\cdot)\in
L^{2}({\mathbb{S}^{2}})\text{ , }\ f(\cdot)\neq 0\text{ .}
\end{equation*}
\end{remark}

\section{Main Results} \label{MainResults}

Throughout this paper, we shall assume to be able to observe the projections
of the fields on the orthonormal basis $\left\{ Y_{\ell m}\right\} ,$ i.e.,
we assume to observe
\begin{equation*}
a_{\ell ,m}(t):=\int_{\mathbb{S}^{2}}T(x,t)Y_{\ell ,m}(x)dx\text{ ,} \qquad %
t=1,\dots,n\text{ .}
\end{equation*}%
The estimator we shall focus on is a form of least square regression on an
increasing subset of the orthonormal system $\left\{ Y_{\ell ,m}\right\} ;$
more precisely, we shall define $k(\cdot ):=(k_{1}(\cdot ),\cdots
,k_{p}(\cdot ))^{\prime }$ for the vector of nuclear kernels, a growing
sequence of integers $L_{N},$ $L_{N}\rightarrow \infty $ as $N\rightarrow
\infty ;$ and a vector of estimators%
\begin{equation}  \label{kernel_estimator}
\widehat{k}_{N}(\cdot ):=(\widehat{k}_{1;N}(\cdot ),\dots ,\widehat{k}%
_{p;N}(\cdot ))^{\prime }=\underset{{k}(\cdot )\,\in \,\mathcal{P}_{N}^{p}}{%
\arg \min }\sum_{t=1}^{N}\left\Vert T_{t+p}-\sum_{j=1}^{p}\Phi
_{j}T_{t+p-j}\right\Vert _{L^{2}(\mathbb{S}^{2})}^{2}\text{ },
\end{equation}%
where $N:= n - p$, $N>p$, and $\mathcal{P}_{N}^{p}$ is the Cartesian product of $p$
copies of
\begin{equation*}
\mathcal{P}_{N}=span\left\{ \frac{2\ell +1}{4\pi }P_{\ell }(\cdot ),\ \ell
\leq L_{N}\right\} \text{ .}
\end{equation*}%
As common in the autoregressive context, we drop the first $p$ observations
when computing our estimators, in order to avoid initialization issues. We
shall write $L_{N}(\cdot)$ for the function $L_{N}(\cdot ):[-1,1]\rightarrow
\mathbb{R}$,
\begin{equation}
L_{N}(z)=\sum_{\ell =0}^{L_{N}}\frac{2\ell +1}{16\pi ^{2}}P_{\ell }^{2}(z)%
\text{ },\qquad z\in \lbrack -1,1]\text{ }.  \label{L_N(z)}
\end{equation}
Note that%
\begin{equation*}
L_{N}(1)=L_{N}(-1)=\sum_{\ell =0}^{L_{N}}\frac{2\ell +1}{16\pi ^{2}}=\frac{%
(L_{N}+1)^{2}}{16\pi ^{2}}\text{ ;}
\end{equation*}
on the other hand, for $z\in (-1,1)$ we have the identity (see \cite%
{Szego,GR})%
\begin{equation*}
\sum_{\ell =0}^{L_{N}}\frac{2\ell +1}{16\pi ^{2}}P_{\ell }^{2}(z)=\frac{%
L_{N}+1}{16\pi ^{2}}\left[ P_{L_{N}+1}^{\prime
}(z)P_{L_{N}}(z)-P_{L_{N}}^{\prime }(z)P_{L_{N}+1}(z)\right] \text{ ;}
\end{equation*}%
it is then possible to show that (see Lemma 4 in the Supplementary Material) 
\begin{equation}
L_{N}(z) \simeq \frac{2L_{N}}{\pi\sqrt{1-z^{2}}}\ , \qquad \text{as $L_{N}
\rightarrow \infty$ ,}  \label{lower L}
\end{equation}
where $\simeq$ indicates that the ratio of left- and right-hand sides
converges to unity.

For our results to follow, we need slightly stronger assumptions on the
"high frequency" behaviour of the kernels $k_{j}(\cdot ).$ More precisely,
we shall introduce the following:

\begin{condition}[Smoothness]
\label{smoothness}For all $j=1,\dots,p$ there exists positive
constants $\beta _{j},\gamma _{j}$ such that%
\begin{equation}
|\phi _{\ell;j}|\leq \frac{\gamma _{j}}{\ell ^{\beta _{j}}}\ ,\qquad \beta
_{j}>1\ ,\ \ell >0\text{ . }  \label{decay::beta_j}
\end{equation}
We let $\beta _{\ast }=\min_{j\in \{1,\dots ,p\}}\beta _{j}$. We shall say
that this condition is satisfied in the \emph{strong} sense if $\beta_j >2,
\ j=1,\dots,p$.
\end{condition}

\begin{remark}
It is readily seen that Condition \ref{smoothness} leads to Hilbert-Schmidt
operators, since it implies $\sum_{\ell }(2\ell +1)\phi _{\ell
;j}^{2}<\infty $, $j=1,\dots,p$; whereas the strong version Condition \ref%
{smoothness} is specific for nuclear operators, since it entails $\sum_{\ell
}(2\ell +1)|\phi _{\ell ;j}|<\infty $, $j=1,\dots,p$, see again \cite{Hsing}.
\end{remark}

\begin{remark}
Condition \ref{smoothness} is easily interpretable in terms of the
regularity of each kernel $k_{j}(\cdot )$. Indeed, in \cite{LangSchwab} it
is shown that
\begin{equation*}
\sum_{\ell =0}^{\infty }|\phi _{\ell;j}|^{2}\frac{2\ell +1}{4\pi }(1+\ell
^{2\eta })<\infty
\end{equation*}%
implies integrability of the first $\eta $ derivatives of $k_{j}(\cdot ),$
i.e., $k_{j}(\cdot )$ belongs to the Sobolev space $W_{1,\eta }$.
\end{remark}

Our first result refers to the asymptotic consistency of the kernel
estimators that we just introduced.

\begin{theorem}[Consistency]\label{consistency}
Consider $\widehat{{k}}_{N}(\cdot )$ in Equation \eqref{kernel_estimator}.
Under Conditions \ref{stationarity}, \ref{Identifiability} and \ref%
{smoothness}, for $L_{N}\sim N^{d},\ 0<d<1$, we have that
\begin{equation}
{\mathbb{E}}\left[ \int_{-1}^{1}\left\Vert \widehat{{k}}_{N}(z)-{k}%
(z)\right\Vert ^{2}dz\right] =\mathcal{O}\left( N^{d-1}+N^{2d(1-\beta _{\ast
})}\right) \text{ .}  \label{imse}
\end{equation}%
Moreover, under Conditions \ref{stationarity}, \ref{Identifiability} and \ref%
{smoothness} (in the \emph{strong} sense), for $L_{N}\sim N^{d},\ 0<d<\frac{1%
}{3}$,%
\begin{equation*}
\mathbb{E}\left[ \sup_{z\in \lbrack -1,1]}\left\Vert \widehat{{k}}_{N}(z)-{k}%
(z)\right\Vert \right] =\mathcal{O}\left( N^{(3d-1)/2}+N^{d(2-\beta _{\ast
})}\right) \text{ .}
\end{equation*}
\end{theorem}

\begin{remark}[Optimal choice of $d$]\label{optimal_choice}
The optimal choice of $d,$ in
terms of the best convergence rates, is given by $d^{\ast }=\frac{1}{2\beta
_{\ast }-1},$ leading to the exponents $\frac{2-2\beta _{\ast }}{2\beta
_{\ast }-1}$ and $\frac{2-\beta _{\ast }}{2\beta _{\ast }-1},$ respectively$.
$ Heuristically, the result can be explained as follows: larger values of $%
\beta _{\ast }$ entail higher regularity/smoothness properties of the
kernels to be estimated; as usual in nonparametric estimation, more regular
functions can be estimated with better convergence rates, as the bias term
is controlled more efficiently. Indeed, for $d=d^{\ast }$ and $\beta _{\ast
}\rightarrow \infty $, the mean square error approximates the parametric
rate $1/N$, as expected.
\end{remark}

\begin{remark}[Plug-in estimates] For applications to empirical data, the optimal rate can
be implemented by means of plug-in techniques, i.e., estimating (under
additional regularity conditions) the value of the parameter $\beta _{\ast }$
by means of first step-estimators of the coefficients $\phi _{\ell ,j}.$ Let
us sketch the main ideas for this approach, omitting some details for
brevity. Consider for simplicity the $SPHAR(1)$ case, and let us make
Condition \ref{smoothness} stronger by assuming that%
\begin{equation*}
|\phi _{\ell }|=\frac{\gamma }{\ell ^{\beta }}+o\left (\frac{1}{\ell ^{\beta }} \right)%
\text{ },\qquad \text{some }\gamma >0\text{ , }\beta >1\ ,\ \forall %
\ell >0\text{ . }
\end{equation*}%
Consider the estimator%
\begin{equation*}
\widehat{\phi }_{\ell ,N}:=\frac{\sum_{t}a_{\ell ,m}(t-1)a_{\ell ,m}(t)}{%
\sum_{t}a_{\ell ,m}^{2}(t-1)}\text{ ,} \qquad \ell =0,1,2,\dots,
\end{equation*}%
from which we can now build the pseudo log-regression model%
\begin{align*}
\log \widehat{\phi }_{\ell ,N}^{2}=&\log \frac{\widehat{\phi }_{\ell ,N}^{2}}{%
\gamma ^{2}\ell ^{-2\beta }}+\log \left( \gamma ^{2}\ell ^{-2\beta
}\right) =\log \left ( \gamma ^{2}\right) -2\beta \log \ell
+v_{\ell }\\
v_{\ell }:=&\log \frac{\widehat{\phi }_{\ell ,N}^{2}}{\gamma ^{2}\ell ^{-2\beta }}
\text{ ,} \qquad \ell =0,1,2,\dots,
\end{align*}%
where the "regression residuals" $\left\{ v_{\ell }\right\}$
are independent over $\ell$, with asymptotically mean zero and bounded
variance as $N\rightarrow \infty .$  It is then possible to study the
asymptotic consistency of the OLS-like estimator (see also \cite%
{Robinson1995} for the related log-periodogram estimator)%
\begin{equation*}
\widehat{\beta }_{N}:=-\frac{\sum_{\ell }\left\{ \log \ell
\times \log \widehat{\phi }_{\ell ,N}^{2}\right\} }{2\sum_{\ell }\left\{
\log \ell \right\} ^{2}}\text{ .}
\end{equation*}%
The optimal rates can then be consistently estimated by means of the plug-in
estimates $\widehat{d}_{N}^{\ast }=\frac{1}{2\widehat{\beta }_{N}-1}.$

A more rigorous and complete investigation on these issues is currently in preparation and is not reported here
for brevity's sake.
\end{remark}

Our second result refers to a Quantitative Central Limit Theorem for the
kernel estimators. Consider $\widehat{{k}}_{N}(\cdot )$ in Equation %
\eqref{kernel_estimator} and, for any $m\in \mathbb{N}$, any $z_{1},\dots
,z_{m}\in (-1,1)$, $z_{1}\neq \cdots \neq z_{m}$, define the $mp\times 1$
vectors
\begin{equation*}
{K}_{N}={K}_{N}(z_{1},z_{2},\dots,z_{m}):=%
\begin{pmatrix}
\sqrt{\frac{N}{L_{N}(z_{1})}}\left( \widehat{{k}}_{N}(z_{1})-{k}%
(z_{1})\right) \\
\vdots \\
\sqrt{\frac{N}{L_{N}(z_{m})}}\left( \widehat{{k}}_{N}(z_{m})-{k}%
(z_{m})\right)%
\end{pmatrix}%
\text{ ,}\quad {Z}\, \overset{d}{{=}}\, \mathcal{N}_{mp}\left( {0}%
_{mp},I_{mp}\right) \text{ .}
\end{equation*}

\begin{theorem}
\label{main2} Under Conditions \ref{stationarity}, \ref{Identifiability} and %
\ref{smoothness} (in the \emph{strong} sense), for $L_{N}\sim N^{d},\ d>%
\frac{1}{2\beta _{\ast }-2}$, we have that
\begin{equation*}
d_{W}({Z},{K}_{N})=\mathcal{O}\left( N^{-1/2}+N^{1/2+d(1-\beta _{\ast } )}+N^{-d}\log
N\right) \text{ .}
\end{equation*}
\end{theorem}

\begin{remark}
  It is easy to see that the bound in Theorem \ref{main2} can also be
expressed as%
 \begin{equation*}
d_{W}({Z},{K}_{N})=\mathcal{O}\left( N^{-1/2}+N^{-d_{\ast }}\log N\right)
\text{ ,} \quad d_{\ast }:=\frac{1}{2\beta _{\ast }-4}\text{ .}
\end{equation*}
\end{remark}

An immediate Corollary is the following.

\begin{corollary}
Under the same Conditions and notation as in Theorem \ref{main2}, for any
fixed $z\in \lbrack -1,1]$, we have that
\begin{equation*}
\sqrt{\frac{N}{L_{N}(z)}}\left( \widehat{{k}}_{N}(z)-{k}(z)\right)
\rightarrow \mathcal{N}_{p}\left( {0}_{p},I_{p}\right) ,\qquad N\rightarrow
\infty .
\end{equation*}
\end{corollary}

\begin{remark}
Plug-in procedures can be exploited to determine the choice of the "bandwidth" parameter $d$ which yields the optimal convergence rate in Wasserstein distance. As usual, the values of $d$ that guarantee asymptotic normality do not minimize the mean squared error; in fact, we have that $d^{\ast }=\frac{1}{%
2\beta _{\ast }-1}<\frac{1}{2\beta _{\ast }-2},$ which is the minimal value
of $d$ for Theorem \ref{main2} to hold. Indeed, asymptotic Gaussianity
requires undersmoothing, i.e., a value of $d$ which makes the asymptotic
bias negligible, rather than of the same order as the variance. Once again
the rate can be taken to approach $N^{-1/2}$ for $\beta _{\ast }\rightarrow
\infty .$
\end{remark}

For our third and final result, we need to strengthen the conditions on the
regularity of the autoregressive kernels.

\begin{condition}
\label{FiniteExpansion} The kernel $k_{j}(\cdot)$ admits a final expansion in
the Legendre basis, i.e., there exist an (arbitrary large but finite)
integer $L>0$ such that%
\begin{equation*}
\int_{\mathbb{S}^{2}}k_{j}(x)P_{\ell }(x)dx=0 \ , \qquad \text{for all } j=1,\dots,p%
\text{ and }\ell >L\text{ .}
\end{equation*}
\end{condition}

Condition \ref{FiniteExpansion} clearly implies that there exist finite
integers $L_{1},\dots,L_{p}\leq L$ such that%
\begin{equation*}
k_{j}(z)=\sum_{\ell =0}^{L_{j}}\frac{2\ell +1}{4\pi }\phi _{\ell ;j}P_{\ell
}(z)\text{ ,} \qquad z \in [ -1,1]\ , \ j=1,\dots,p\text{ ;}
\end{equation*}%
we also need to introduce, for $\ell=0,1,2,\dots,$ the $p\times p$ autocovariance matrix
\begin{equation*}
\Gamma _{\ell }:=
\begin{pmatrix}
C_{\ell} & C_{\ell}(1) & \cdots & C_{\ell}(p-1) \\
C_{\ell }(1) & C_{\ell} & \cdots & C_{\ell}(p-2) \\
\vdots & \vdots & \ddots & \vdots \\
C_{\ell} (p-1) & C_{\ell}(p-2) & \cdots & C_{\ell}%
\end{pmatrix} \text{ ,}
\end{equation*}%
and we shall write $W_{p}(\cdot )$ for the zero-mean, $p$-dimensional
Gaussian process with covariance function
\begin{equation*}
\Gamma _{{k}}(z,z^{\prime })=\sum_{\ell =0}^{L}C_{\ell ;Z}\Gamma_{\ell }^{-1}%
\frac{2\ell +1}{16\pi ^{2}}P_{\ell }(z)P_{\ell }(z^{\prime })\text{ .}
\end{equation*}%
We are now able to state our last Theorem.

\begin{theorem}
\label{weakconvergence}Under Conditions \ref{stationarity}, \ref{Identifiability} and \ref{FiniteExpansion}, we have that
\begin{equation*}
\sqrt{N}\left( \widehat{{k}}_{N}(\cdot )-{k}(\cdot )\right) \Longrightarrow
\mathnormal{W}_{p}\left( .\right) \text{ },\qquad N\rightarrow \infty \text{
},
\end{equation*}%
where $\Longrightarrow $ denotes weak convergence in $C_{p}[-1,1]$ (the
space of continuous functions from $[-1,1]$ to $\mathbb{R}^{p},$ with the
standard uniform metric).
\end{theorem}

\begin{remark}
At first sight, it may look surprising that the weak convergence for the
estimators in Theorem \ref{weakconvergence} occurs at a faster rate $\sqrt{N}
$ than the convergence in finite-dimensional distributions of Theorem \ref%
{main2}. This comparison, however, is misleading; indeed, in Theorem \ref%
{main2} we are not assuming the expansion of the kernels to be finite, and
therefore we need to include a growing number of multipoles $L_{N},$ to
ensure that bias terms are asymptotically negligible. On the other hand,
note that weak convergence cannot hold under the conditions of Theorem \ref%
{main2}, as the limiting finite dimensional distributions correspond to
Gaussian independent random variables for any choice of fixed points $(z_{1},\dots,z_{m}):$ no Gaussian process with measurable trajectories can have
these finite-dimensional distributions. The limiting distribution is
characterized by the nuisance parameters $\left\{ C_{\ell },C_{\ell
}(1),\dots,C_{\ell }(p-1),C_{\ell ;Z}\right\};$ for brevity's
sake, estimation of these parameters is deferred to future work.
\end{remark}

\section{Proofs of the Main Results} \label{Proofs}

We now present the main arguments of our proofs, which are based on a number
of technical results collected in the Appendix (Supplementary Material). For $\ell
=0,1,2,\dots,$ it is convenient to introduce the $N(2\ell +1)$-dimensional
vectors
\begin{align*}
&{Y}_{\ell ;N}:=(a_{\ell ,-\ell }(p+1),\dots ,a_{\ell ,\ell }(p+1),\dots
,a_{\ell ,\ell }(n))^{\prime }\text{ , } \\
&\boldsymbol{\varepsilon }_{\ell ;N}:=(a_{\ell ,-\ell ;Z}(p+1),\dots
,a_{\ell ,\ell ;Z}(p+1),\dots ,a_{\ell ,\ell ;Z}(n))^{\prime }\text{ ; }
\end{align*}%
moreover, let us consider the $N(2\ell +1)\times p$ matrix%
\begin{equation*}
X_{\ell ;N}:=\left\{ {Y}_{\ell ;N-1}:{Y}_{\ell ;N-2}:\cdots:{Y}_{\ell
;N-p}\right\} \text{ ,}
\end{equation*}
where%
\begin{equation*}
{Y}_{\ell ;N-j}:=(a_{\ell ,-\ell }(p+1-j),\dots ,a_{\ell ,\ell
}(p+1-j),\dots ,a_{\ell ,\ell }(n-j))^{\prime }\text{ ,} \quad j=1,\dots,p\text{ .}
\end{equation*}

We start from the proof of the consistency results.

\begin{proof}[Proof (Theorem \ref{consistency})]
It is easy to see that we have
\begin{align*}
\widehat{{k}}_{N}(\cdot)
&= \underset{ {k}(\cdot)\, \in\, \mathcal{P}^p_{N}}{ \arg \min}\sum_{t=p+1}^{n}\left\Vert T_{t}-\sum_{j=1}^{p}\Phi
_{j}T_{t-j}\right\Vert _{L^{2}(\mathbb{S}^{2})}^{2} \\
&=\sum_{\ell =0}^{L_{N}}\widehat{\boldsymbol{\phi }}_{\ell ;N}\frac{2\ell +1%
}{4\pi }P_{\ell }(\cdot)\text{ ,}
\end{align*}%
where%
\begin{align*}
\widehat{\boldsymbol{\phi }}_{\ell ;N} &:= ( \widehat{\phi }_{\ell
;N}(1),\dots ,\widehat{\phi }_{\ell ;N}(p) ) ^{\prime } \\
&=\underset{\boldsymbol{\phi }_{\ell } \in \mathbb{R}^p}{\arg \min} \sum_{t=p+1}^{n} \sum_{m=-\ell}^\ell \left (a_{\ell ,m}(t)-\sum_{j=1}^p\phi
_{\ell; j}a_{\ell ,m}(t-j) \right)^{2} \ .
\end{align*}%

Now, let ${r}_N(z)$ be the difference between the kernel and its truncated version
\[
{k}_N(\cdot) = \sum_{\ell=0}^{L_N}\boldsymbol{\phi}_\ell  \frac{2\ell +1}{4\pi} P_\ell(z),
\]
i.e.,
\[
{r}_N(z) = {k}(z) - {k}_N(z)  = \sum_{\ell=L_N+1}^\infty \boldsymbol{\phi}_\ell  \frac{2\ell +1}{4\pi} P_\ell(z)\ ,
\]
where the equality holds in the $L^2$ sense. Then,
\begin{equation}\label{biasvariance}
\Ex \left [ \int_{-1}^1 \left \|\widehat{{k}}_{N}(z) - {k}(z) \right \|^2d z \right] = \Ex \left [ \int_{-1}^1 \left \|\widehat{{k}}_{N}(z) - {k}_N(z) \right \|^2d z \right] +  \int_{-1}^1 \left \|{r}_N(z) \right \|^2d z\ ,
\end{equation}
since $\Ex\left [ \int_{-1}^1 \left \langle \widehat{{k}}_{N}(z) - {k}_N(z), {r}_N(z) \right \rangle d z \right]=0$, from orthogonality of Legendre polynomials.

Now notice that
\begin{align*}
 \int_{-1}^1 \left \|\widehat{{k}}_{N}(z) - {k}_N(z) \right \|^2d z &=   \sum_{\ell =0}^{L_N} \sum_{\ell' =0}^{L_N} \left \langle \widehat{\boldsymbol{\phi}}_{\ell;N} - \boldsymbol{\phi}_\ell,  \widehat{\boldsymbol{\phi}}_{\ell';N} - \boldsymbol{\phi}_{\ell'}\right \rangle \frac{2\ell +1}{4\pi}  \frac{2\ell' +1}{4\pi}  \int_{-1}^1P_\ell (z)  P_{\ell'} (z) d z \nonumber \\
&=  \sum_{\ell =0}^{L_N} \sum_{\ell' =0}^{L_N}  \left \langle \widehat{\boldsymbol{\phi}}_{\ell;N} - \boldsymbol{\phi}_\ell,  \widehat{\boldsymbol{\phi}}_{\ell';N} - \boldsymbol{\phi}_{\ell'}\right \rangle \frac{2\ell +1}{4\pi}  \frac{2\ell' +1}{4\pi} \frac{2}{2\ell + 1} \delta_{\ell}^{\ell'} \nonumber \\
&= \sum_{\ell = 0}^{L_N} \left \| \widehat{\boldsymbol{\phi}}_{\ell;N} - \boldsymbol{\phi}_\ell\right\|^2 \frac{2\ell +1}{8\pi^2}  \ .
\end{align*}
Then, from Lemma 2 in the Supplementary material,  
\begin{align*}
\Ex \left [ \int_{-1}^1 \left \|\widehat{{k}}_{N}(z) - {k}_N(z) \right \|^2d z \right]
&=  \sum_{\ell = 0}^{L_N} \Ex \left \| \widehat{\boldsymbol{\phi}}_{\ell;N} - \boldsymbol{\phi}_\ell\right\|^2 \frac{2\ell +1}{8\pi^2}  \le const \frac{L_N + 1}{N}\ .
\end{align*}

On the other hand,
\begin{align*}
\int_{-1}^1 \left \|{r}_N(z) \right \|^2d z &=\sum_{\ell = L_N+1}^\infty  \sum_{\ell' = L_N+1}^\infty  \left \langle  \boldsymbol{\phi}_\ell, \boldsymbol{\phi}_{\ell'}\right \rangle \frac{2\ell +1}{4\pi}    \frac{2\ell' +1}{4\pi}  \int_{-1}^1P_\ell (z)  P_{\ell'} (z) d z \nonumber \\
&= \sum_{\ell = L_N+1}^\infty \sum_{\ell' = L_N+1}^\infty  \left \langle  \boldsymbol{\phi}_\ell, \boldsymbol{\phi}_{\ell'}\right \rangle  \frac{2\ell +1}{4\pi}  \frac{2\ell' +1}{4\pi} \frac{2}{2\ell + 1} \delta_{\ell}^{\ell'} \nonumber \\
&= \sum_{\ell = L_N+1}^\infty  \| \boldsymbol{\phi}_\ell \|^2 \frac{2\ell +1}{8 \pi^2}.
\end{align*}
Therefore, under Condition \ref{smoothness} and for $L_N\sim N^d, \ 0<d<1$, we have
\begin{equation*}
\int_{-1}^1 \left \|{r}_N(z) \right \|^2d z  = \mathcal{O} \left ( N^{2d(1-\beta_{\ast})} \right),
\end{equation*}
and
\begin{equation*}
\Ex \left [ \int_{-1}^1 \left \|\widehat{{k}}_{N}(z) - {k}(z) \right \|^2d z \right] = \mathcal{O} \left (  N^{d-1} + N^{2d(1-\beta_{\ast})}    \right ),
\end{equation*}
where $\beta_{\ast} = \min_{j \in \{1,\dots,p\}} \beta_j$, as claimed.

Under the strong version of Condition \ref{smoothness}, each kernel $k_j(\cdot)$ is defined for all $z \in [-1,1]$ as the pointwise limit of its expansion in terms of Legendre polynomials and
\begin{equation*}
\Ex \left [  \sup_{z \in [-1,1]} \left \| \widehat{{k}}_{N}(z) - {k}(z)\right \|   \right ] \le \Ex \left [   \sup_{z \in [-1,1]} \left \| \widehat{{k}}_{N}(z) - {k}_N(z)\right \| \right ]   +   \sup_{z \in [-1,1]} \left \|{r}_N(z)\right \| \ ,
\end{equation*}
by the triangle inequality. Hence, for the first component we have
\begin{align*}
\Ex \, \left [ \sup_{z \in [-1,1]} \left \| \sum_{\ell =0}^{L_N} \left ( \widehat{\boldsymbol{\phi}}_{\ell;N} - \boldsymbol{\phi}_\ell\right) \frac{2\ell +1}{4\pi} P_\ell (z) \right \| \right ] &\le \sum_{\ell =0}^{L_N}  \Ex \left \|  \widehat{\boldsymbol{\phi}}_{\ell;N} - \boldsymbol{\phi}_\ell  \right \|\frac{2\ell +1}{4\pi} \nonumber \\
&\le const  \sum_{\ell =0}^{L_N}\frac{\sqrt{2\ell +1}}{\sqrt{N}} \nonumber\\
&\le const \frac{(L_N+1)^{3/2}}{\sqrt{N}} \ ,
\end{align*}
again in view of Lemma 2 in the Appendix (Supplementary Material)
and the Cauchy-Schwartz inequality. On the other hand
\begin{align*}
 \sup_{z \in [-1,1]} \left \| {r}_N(z)  \right \| &\le \sum_{\ell=L_N+1}^\infty \| \boldsymbol{\phi}_\ell \| \frac{2\ell +1}{4\pi} \ .
\end{align*}
Therefore, again under the strong version of Condition \ref{smoothness} and for $L_N\sim N^d, \ 0<d<\frac13$, we have
\begin{equation*}
 \sup_{z \in [-1,1]} \left \| {r}_N(z)  \right \| = \mathcal{O} \left ( N^{d(2-\beta_{\ast})}  \right)
\end{equation*}
and thus
\begin{equation*}
\mathbb{E} \left [ \sup_{z \in [-1,1]} \left \|\widehat{{k}}_{N}(z) -
{k}(z) \right \| \right] = \mathcal{O} \left ( N^{(3d-1)/2} +
N^{d(2-\beta_{\ast})} \right )\text{ .}
\end{equation*}
as claimed.
\end{proof}

We are now in the position to establish the Quantitative Central Limit
Theorem.

\begin{proof}[Proof (Theorem \ref{main2})] Let us recall that the minimizing estimator takes the form
\begin{align*}
\widehat{{k}}_{N}(\cdot)
&= \underset{ {k}(\cdot)\, \in\, \mathcal{P}^p_{N}}{ \arg \min}\sum_{t=p+1}^{n}\left\Vert T_{t}-\sum_{j=1}^{p}\Phi
_{j}T_{t-j}\right\Vert _{L^{2}(\mathbb{S}^{2})}^{2}\\
&=\sum_{\ell =0}^{L_{N}}\widehat{\boldsymbol{\phi }}_{\ell ;N}\frac{2\ell +1%
}{4\pi }P_{\ell }(\cdot)\text{ ,}
\end{align*}%
where%
\begin{align*}
\widehat{\boldsymbol{\phi }}_{\ell ;N}&=\underset{\boldsymbol{\phi }_{\ell } \in \mathbb{R}^p}{\arg \min} \sum_{t=p+1}^{n} \sum_{m=-\ell}^\ell \left (a_{\ell ,m}(t)-\sum_{j=1}^p\phi_{\ell ;j}a_{\ell ,m}(t-j) \right)^{2}\\
&= (X_{\ell ;N}^{\prime }X_{\ell ;N})^{-1}X_{\ell ;N}^{\prime }{Y}%
_{\ell ;N}\text{ }=\boldsymbol{\phi }_{\ell }+(X_{\ell ;N}^{\prime }X_{\ell
;N})^{-1}X_{\ell ;N}^{\prime }\boldsymbol{\varepsilon }_{\ell ;N}\text{ }.\text{
} 
\end{align*}%
We shall introduce some more notation: 
\begin{equation*}
A_{\ell ;N}:=\frac{1}{C_{\ell }N(2\ell +1)}X_{\ell ;N}^{\prime }X_{\ell ;N}\ ,
\qquad   \Sigma_\ell := \Ex [ A_{\ell ;N} ] = \frac{\Gamma_\ell}{C_\ell}\ ,
\end{equation*}
and
\begin{equation*}
{B}_{\ell ;N}:=\frac{1}{C_{\ell }\sqrt{N(2\ell +1)}}X_{\ell
;N}^{\prime }\boldsymbol{\varepsilon }_{\ell ;N}\ . 
\end{equation*}
Therefore
\begin{equation*}\label{pierpaolo}
\sqrt{N(2\ell+1)} \left (\widehat{\boldsymbol{\phi }}_{\ell ;N}-\boldsymbol{\phi }_{\ell } \right ) = A_{\ell ;N}^{-1}{B}_{\ell ;N} \ .
\end{equation*}
Heuristically, the proof of the Quantitative Central Limit Theorem can be described as follows: in order to be able to exploit Stein-Malliavin techniques, we need to deal with variables belonging to some $q$-th order chaos; now the ratio above does not fulfill this requirement, because $A_{\ell ;N}^{-1}$ is a random quantity which does not belong to any $\mathcal{H}_q$. On the other hand, componentwise we have ${B}_{\ell ;N} \in \mathcal{H}_2$, for each $\ell$. We shall then show that it is possible to replace $A_{\ell ;N}^{-1}$ by its (deterministic) probability limit $\Sigma_{\ell}^{-1}$, without affecting asymptotic results; because our kernel estimators will be written as linear combinations of $\widehat{\boldsymbol{\phi }}_{\ell ;N}$, the proof can be completed by a careful investigation of multivariate fourth-order cumulants.

Let us now make the previous argument rigorous. Let ${K}_{N}$ and ${U}_{N}$ be two $mp$-dimensional random
vectors, defined as
\begin{equation*}
{K}_{N}:=
\begin{pmatrix}
\sqrt{\frac{N}{L_{N}(z_{1})}}\left( \widehat{{k}}_{N}(z_{1})-{k%
}(z_{1})\right)  \\
\vdots  \\
\sqrt{\frac{N}{L_{N}(z_{m})}}\left( \widehat{{k}}_{N}(z_{m})-{k%
}(z_{m})\right)
\end{pmatrix}%
,
\end{equation*}%
and
\begin{equation*}
{U}_{N}=%
\begin{pmatrix}
{U}_{N}(z_{1}) \\
\vdots  \\
{U}_{N}(z_{m})%
\end{pmatrix}%
:=%
\begin{pmatrix}
\frac{1}{\sqrt{L_{N}(z)}}\sum_{\ell =0}^{L_{N}}\Sigma _{\ell }^{-1}{B}%
_{\ell ;N}\frac{\sqrt{2\ell +1}}{4\pi }P_{\ell }(z_{1}) \\
\vdots  \\
\frac{1}{\sqrt{L_{N}(z_m)}}\sum_{\ell =0}^{L_{N}}\Sigma _{\ell }^{-1}{B}%
_{\ell ;N}\frac{\sqrt{2\ell +1}}{4\pi }P_{\ell }(z_{m})%
\end{pmatrix}%
.
\end{equation*}%
In particular, ${\mathbb{E}}[{U_{N}}]={0}_{mp}$ and ${\mathbb{E%
}}[{U}_{N}{U}_{N}^{\prime }]=V_{N}$, where $V_{N}$ is a block
matrix whose generic $ij$-th block, $i,j\in \{1,\dots ,m\}$, is given by
\begin{align*}
& V_N(i,j)={\mathbb{E}}[{U}_{N}(z_{i}){U}_{N}^{\prime }(z_{j})]\\
=&\frac{1}{\sqrt{L_{N}(z_{i})}}\frac{1}{\sqrt{L_{N}(z_{j})}}\sum_{\ell =0}^{L_{N}}\frac{%
C_{\ell; Z}}{C_{\ell }}\Sigma _{\ell }^{-1}\frac{2\ell +1}{16\pi ^{2}}%
P_{\ell }(z_{i})P_{\ell }(z_{j})\ .
\end{align*}

Now, consider ${Z}\overset{d}{=}\mathcal{N}_{mp}({0}_{mp},I_{mp})$ and $%
{Z}_{N}\overset{d}{=} \mathcal{N}_{mp}({0}_{mp},V_{N})$. Applying
the triangle inequality twice, it follows that
\begin{align*}
d_{W}({Z},{K}_{N})& \leq d_{W}({Z},{U}%
_{N})+d_{W}({U}_{N},{K}_{N})  \notag \\
& \leq d_{W}({Z},{Z}_{N})+d_{W}({Z}_{N},{U}%
_{N})+d_{W}({U}_{N},{K}_{N})\ .
\end{align*}%

We recall from \cite{noupebook}, p. 126, Equation 6.4.2 that
\begin{equation*}
d_{W}({Z},{Z}_{N})\leq \sqrt{mp}\min \{\Vert V_{N}^{-1}\Vert _{%
\text{op}}\Vert V_{N}\Vert _{\text{op}}^{1/2},1\}\Vert V_{N}-I_{mp}\Vert _{%
\text{HS}}\ ,
\end{equation*}%
where $\Vert A\Vert _{\text{HS}}=\sqrt{\text{Tr}(A^{\prime }A)}$, and we observe that
\begin{equation}\label{VN::hs}
\Vert V_{N}-I_{mp}\Vert _{\text{HS}} \leq mp\Vert V_{N}-I_{mp}\Vert
_{\infty }  = \mathcal{O}\left ( N^{-d} \log N\right) \ ,
\end{equation}%
from Lemmas 3 and 4
in the Supplementary Material. Indeed, for every $i \in \{1,\dots, m\}$,
\begin{align*}
\| V_N(i,i) - I_p \|_{\text{HS}}& \leq \frac{const}{L_{N}+1}\sum_{\ell =0}^{L_{N}}\left \Vert\frac{C_{\ell;Z}}{C_\ell} \Sigma _{\ell
}^{-1}-I_{p} \right \Vert _{\infty }(2\ell +1)  \notag \\
& \leq \frac{const}{L_{N}+1} \ ;
\end{align*}
the logarithmic term comes from Equation (8)
in the Supplementary Lemma 4.
Equation \eqref{VN::hs} entails that $V_{N}\rightarrow I_{mp}$, thus we have
$\Vert V_{N}^{-1}\Vert _{\text{op}}\Vert V_{N}\Vert _{\text{op}%
}^{1/2}\rightarrow 1$, as $N\rightarrow \infty $, and
\begin{equation}
d_{W}({Z},{Z}_{N})=\mathcal{O}\left( N^{-d} \log N%
\right)\ .   \label{triangle::1}
\end{equation}%

Let us recalll again from \cite{noupebook}, p. 122 (second point of Theorem 6.2.2) that
\begin{equation*}
d_{W}({Z}_{N},{U}_{N})\leq \sqrt{mp}\Vert V_{N}^{-1}\Vert _{%
\text{op}}\Vert V_{N}\Vert _{\text{op}}^{1/2}m({U}_{N})\ ,
\end{equation*}%
where
\begin{align*}
m({U}_{N})& =2mp\sum_{i=1}^{m}\sum_{j=1}^{p}\sqrt{\textrm{Cum}_4 \left [ \frac{1}{\sqrt{L_N(z_i)}}\sum_{\ell=0}^{L_N}  \tilde{b}_{\ell ;N}(j) \frac{\sqrt{2\ell +1}}{4\pi} P_\ell (z_i)\right] } \ ,
\end{align*}%
$\tilde{b}_{\ell ;N}(j)$ being the $j$-th element of $\Sigma_\ell^{-1}{B}_{\ell;N}$. 
Moreover, for the $j$-th element of $\Sigma_\ell^{-1}{B}_{\ell;N}$ we have
\begin{align*}
\text{Cum}_4\left [\tilde{b}_{\ell ;N}(j) \right]= \frac{6}{N(2\ell+1)} \left (\frac{C_{\ell; Z}}{C_\ell} s_\ell(j,j) \right)^2 \ ,
\end{align*}
see Equation (4)
in Lemma 1.
In addition,
\begin{align*}
&\textrm{Cum}_4 \left [ \frac{1}{\sqrt{L_N(z_i)}}\sum_{\ell=0}^{L_N}  \tilde{b}_{\ell ;N}(j) \frac{\sqrt{2\ell +1}}{4\pi} P_\ell (z_i)\right]\\
=& \frac{1}{L_N^2(z_i)} \sum_{\ell=0}^{L_N}  \textrm{Cum}_4\left [\tilde{b}_{\ell ;N}(j)\right ] \frac{(2\ell +1)^2}{(4\pi)^4} P^4_\ell (z_i)\ ,
\end{align*}
in view of the independence across different multipoles $\ell$. Therefore
\begin{align*}
&\textrm{Cum}_4 \left [ \frac{1}{\sqrt{L_N(z_i)}}\sum_{\ell=0}^{L_N}\tilde{b}_{\ell ;N}(j) \frac{\sqrt{2\ell +1}}{4\pi} P_\ell (z_i)\right]	\\
=& \frac{6}{NL_N^2(z_i)} \sum_{\ell=0}^{L_N}    \left (\frac{C_{\ell; Z}}{C_\ell}s_\ell(j,j) \right)^2  \frac{2\ell +1}{(4\pi)^4} P^4_\ell (z_i) \nonumber \\
\le& \frac{6}{NL_N^2(z_i)} \sum_{\ell=0}^{L_N}    \left  [\frac{C_{\ell; Z}}{C_\ell} \tr(\Sigma_\ell^{-1})\right ]^2 \frac{2\ell +1}{(4\pi)^4} P^4_\ell (z_i) \nonumber \\
\le& \frac{const}{N(L_N+1)^2} \sum_{\ell=0}^{L_N}  (2\ell + 1)P^4_\ell(z_i) \ .
\end{align*}
Thus, we have%
\begin{equation*}
m({U}_{N}) \leq const\frac{m^{2}p^{2}}{L_{N}+1} \sqrt{\frac{\log N}{N}} \ ,
\end{equation*}%
and
\begin{equation}
d_{W}({Z}_{N},{U}_{N})=\mathcal{O}\left(N^{-(d+1/2)} (\log N )^{1/2}%
\right) \ .   \label{triangle::2}
\end{equation}%

Now, consider the decomposition
\begin{align*}
\sqrt{\frac{N}{L_N(z)}}\left ( \widehat{{k}}_{N}(z) - {k}(z) \right) &=\frac{1}{\sqrt{L_N(z)}}\sum_{\ell =0}^{L_N} \sqrt{N(2\ell+1)} \left ( \widehat{\boldsymbol{\phi}}_{\ell;N} - \boldsymbol{\phi}_\ell\right) \frac{\sqrt{2\ell +1}}{4\pi} P_\ell (z) \nonumber\\
&- \sqrt{\frac{N}{L_N(z)}}\sum_{\ell=L_N+1}^\infty \boldsymbol{\phi}_\ell  \frac{2\ell +1}{4\pi} P_\ell(z)  \nonumber \\
&= \frac{1}{\sqrt{L_N(z)}}\sum_{\ell =0}^{L_N}  \Sigma^{-1}_\ell{B}_{\ell;N} \frac{\sqrt{2\ell +1}}{4\pi} P_\ell (z) \nonumber \\
&+  \frac{1}{\sqrt{L_N(z)}}\sum_{\ell =0}^{L_N}  [A_{\ell;N}^{-1} - \Sigma_\ell^{-1}]{B}_{\ell;N}  \frac{\sqrt{2\ell +1}}{4\pi} P_\ell (z) \nonumber \\
&-\sqrt{\frac{N}{L_N(z)}}\sum_{\ell=L_N+1}^\infty \boldsymbol{\phi}_\ell  \frac{2\ell +1}{4\pi} P_\ell(z)\ .
\end{align*}
Without loss of generality, we shall focus on the case $m=1$; the more general argument is basically identical, with a slightly more cumbersome notation. For $z \in (-1,1)$,
\begin{align*}
&\left \|   \frac{1}{\sqrt{L_N(z)}}\sum_{\ell =0}^{L_N}  [A_{\ell;N}^{-1} - \Sigma_\ell^{-1}]{B}_{\ell;N}   \frac{\sqrt{2\ell +1}}{4\pi} P_\ell (z)  \right \|\\
\le&  \frac{const}{\sqrt{L_N+1}}  \sum_{\ell =0}^{L_N}  \left \|[A_{\ell;N}^{-1} - \Sigma_\ell^{-1}]{B}_{\ell;N}  \right \| \sqrt{2\ell +1} |P_\ell (z)|\ ,
\end{align*}
and then, also using Hilb's equation,
\begin{align}\label{resto::1}
&\Ex \, \left [  \left \|   \frac{1}{\sqrt{L_N(z)}}\sum_{\ell =0}^{L_N}  [A_{\ell;N}^{-1} - \Sigma_\ell^{-1}]{B}_{\ell;N}   \frac{\sqrt{2\ell +1}}{4\pi} P_\ell (z)  \right \| \right ] \nonumber \\
\le&  \frac{const}{\sqrt{L_N+1}}  \sum_{\ell =0}^{L_N}  \Ex \left \|[A_{\ell;N}^{-1} - \Sigma_\ell^{-1}]{B}_{\ell;N}  \right \|  \sqrt{2\ell +1} |P_\ell (z)| \nonumber \\
\le&  \frac{const}{\sqrt{L_N+1}} \sum_{\ell=0}^{L_N} \frac{1}{\sqrt{N(2\ell+1)}} \sqrt{2\ell +1} |P_\ell (z)| \nonumber \\
=&\, \mathcal{O}\left(\frac{1}{\sqrt{N}}\right),
\end{align}
where for the second inequality we have exploited the Supplementary Lemma 2.
Likewise,
\begin{align}\label{resto::2}
 \left \|  \sqrt{\frac{N}{L_N(z)}}\sum_{\ell=L_N+1}^\infty \boldsymbol{\phi}_\ell  \frac{2\ell +1}{4\pi} P_\ell(z)  \right \| &\le const \sqrt{\frac{N}{L_N+1}} \sum_{\ell=L_N+1}^\infty \| \boldsymbol{\phi}_\ell \| (2\ell+1)|P_\ell(z)| \nonumber \\
 &\le const \sqrt{\frac{N}{L_N+1}} \sum_{\ell=L_N+1}^\infty \| \boldsymbol{\phi}_\ell \| \sqrt{2\ell+1} \nonumber \\
&=\mathcal{O}\left( \frac{1}{N^{d(\beta_{*}-1)-1/2}} \right)  \ .
\end{align}
From Equations \eqref{resto::1} and \eqref{resto::2},
\begin{equation}
d_{W}({U}_{N},{K}_{N})=\mathcal{O}\left( N^{-1/2}+N^{1/2+d(1-\beta_{*})} \right) \  .  \label{triangle::3}
\end{equation}%

In the end, combining Equations \eqref{triangle::1}, \eqref{triangle::2} and \eqref{triangle::3},
it holds that
\begin{equation*}
d_{W}({Z},{K}_{N})=\mathcal{O}\left( N^{-1/2}+N^{1/2+d(1-\beta_{*})}\right)  \ .
\end{equation*}%
Note that the constant in this bound may depend on the choice of $m$ and $z_{1},\dots ,z_{m}$.
\end{proof}

We can now give the proof of the third (and final) result.
\begin{proof}[Proof (Theorem \ref{weakconvergence})] We have that, for $z \in [-1,1]$,
\begin{align}\label{weak_kernel_p}
\sqrt{N} \left ( \widehat{{k}}_{L;N}(z) - {k}_L(z) \right) &= \sum_{\ell =0}^L \sqrt{N(2\ell+1)} \left ( \widehat{\boldsymbol{\phi}}_{\ell;N} - \boldsymbol{\phi}_\ell\right) \frac{\sqrt{2\ell +1}}{4\pi} P_\ell (z) \nonumber \\
&= \sum_{\ell =0}^L A_{\ell;N}^{-1}{B}_{\ell;N} \frac{\sqrt{2\ell +1}}{4\pi} P_\ell (z) \nonumber \\
&=\sum_{\ell =0}^L \Sigma_\ell^{-1}{B}_{\ell;N} \frac{\sqrt{2\ell +1}}{4\pi} P_\ell (z) \nonumber \\
&+ \sum_{\ell =0}^L [A_{\ell;N}^{-1} - \Sigma_\ell^{-1}]{B}_{\ell;N}  \frac{\sqrt{2\ell +1}}{4\pi} P_\ell (z)\ .
\end{align}
Then,
\begin{align*}
\sup_{z \in [-1,1]} \left \|  \sum_{\ell =0}^L  [A_{\ell;N}^{-1} - \Sigma_\ell^{-1}]{B}_{\ell;N}  \frac{\sqrt{2\ell +1}}{4\pi} P_\ell (z)  \right \| \le \sum_{\ell =0}^L   \left \|[A_{\ell;N}^{-1} - \Sigma_\ell^{-1}]{B}_{\ell;N} \right \| \frac{\sqrt{2\ell +1}}{4\pi}\ ,
\end{align*}
and hence
\begin{align*}
&\Ex \, \left [ \sup_{z \in [-1,1]} \left \|  \sum_{\ell =0}^L  [A_{\ell;N}^{-1} - \Sigma_\ell^{-1}]{B}_{\ell;N} \frac{\sqrt{2\ell +1}}{4\pi} P_\ell (z)  \right \| \right ] \\
\le& \sum_{\ell =0}^L  \Ex \left \|[A_{\ell;N}^{-1} - \Sigma_\ell^{-1}]{B}_{\ell;N}  \right \| \frac{\sqrt{2\ell +1}}{4\pi} \to 0\ , \qquad N \to \infty \ ,
\end{align*}
in view of the Supplementary Lemma 2.
Then the second part of the sum in \eqref{weak_kernel_p} goes to zero in probability.
Since the sum (over $\ell$) has independent components, we just need to prove that, for each $\ell=0,1,2,\dots,L$, $\{{B}_{\ell;N}P_\ell (\cdot), \ N > 1 \}$ forms a \textit{tight} sequence. Using the tightness criterion given in \cite{Billingsley}, Equation 13.14 on page 143, it is sufficient to show that, for $z_1 \le z \le z_2$,
\begin{align*}
&\Ex \| {B}_{\ell;N}P_\ell (z) - {B}_{\ell;N} P_\ell (z_1)\|  \| {B}_{\ell;N}P_\ell (z_2) - {B}_{\ell;N} P_\ell (z)\|\\
=&\, | P_\ell (z) - P_\ell (z_1)|  | P_\ell (z_2) - P_\ell (z)|  \Ex \|{B}_{\ell;N}\|^2 \nonumber \\
\le&\, p \frac{C_{\ell; Z}}{C_\ell}Q^2_\ell |z - z_1| | z_2 - z | \nonumber \\
\le&\, p \frac{C_{\ell; Z}}{C_\ell} Q^2_\ell (z_2 - z_1)^2\ .
\end{align*}

Convergence of the finite-dimensional distributions is standard and we omit the details, which are closed to those given in the proofs of the previous Theorem. Thus the sequence converges weakly to a zero-mean multivariate Gaussian process with covariance function $\Gamma_{{k}_L}(z,z')=\sum_{\ell=0}^L  C_{\ell; Z}\Gamma_\ell^{-1} \frac{2\ell +1}{16\pi^2}   P_\ell(z) P_\ell (z')$.
\end{proof}

\section{Some Numerical Evidence}  \label{Numerical}

In this section, we present some short numerical results to illustrate the
models and methods that we discussed in this paper.

We stress first that random fields on the sphere cross time can be very
conveniently generated by combining the general features of Python with the
HEALPix software (see \cite{HealPix} and \url{https://healpix.sourceforge.io}). More precisely, HEALPix (which
stands for \emph{Hierarchical Equal Area and iso-Latitude Pixelation}) is a
multi-purpose computer software package for a high resolution numerical analysis of functions on the sphere, based on a clever tessellation scheme: the spherical surface is hierarchically partitioned into curvilinear quadrilaterals of equal area (at a given resolution), distributed on lines of constant latitude, as suggested in the name. In particular, we shall make use
of \texttt{healpy}, which is a Python package based on the HEALPix C++
library. HEALPix was developed to efficiently process Cosmic Microwave
Background data from Cosmology experiments (like \emph{Planck}, \cite{Planck}), but it is now used in many other branches of Astrophysics and applied
sciences.

In short, HEALPix allows to create spherical maps according to the spectral
representation \eqref{spectral_representation}, accepting in input either an
array of random coefficients $\left\{ a_{\ell ,m}\right\} $, or the angular
power spectrum $\left\{ C_{\ell }\right\}$, by means of the routines
\texttt{alm2map} and \texttt{synfast}: in the latter case, random $\left\{
a_{\ell ,m}\right\} $ are generated according to a Gaussian zero
mean distribution with variance $\left\{ C_{\ell }\right\} $. The routine is
extremely efficient and allows to generate maps of resolution up to a few
thousands multipoles in a matter of seconds on a standard laptop computer.

In our case, however, we need random fields where the random spherical harmonics coefficients
have themselves a temporal dependence structure. For this reason, we implemented a simple routine in Python, to simulate Gaussian
$\left\{ a_{\ell ,m}(t)\right\}$ processes, each following an $AR(p)$
dependence structure. These random harmonic coefficients are then uploaded
into HEALPix, to generate maps such as those that are given in Figure 1. In particular, in these two cases we fixed
$L_{\max}=\max (\ell )=30,200$, respectively. Then we generated $\left\{ a_{\ell
,m}(t)\right\}$ according to a stationary $AR(1)$ processes, with
parameters $\phi _{\ell }\simeq const\times \ell ^{-3}$; similarly, we took
here $C_{\ell ;Z}\simeq const\times \ell ^{-2}$. In the figure, we report
the realization for the first 4 periods, simply for illustrative purposes.

We are now in the position to use simulations to validate the previous
results. In our first Tables \ref{Table1}-\ref{Table3}, we report for $B=1000
$ Monte Carlo replications the values of the "variance" and "bias" terms,
i.e., the first and second summand in the mean square equation \eqref{biasvariance}; the second term is actually deterministic, and it is
reported to illustrate the approximation one obtains by cutting the
expansion to a finite multipole value. In the third column, we report, the
actual (squared) $L^{2}$ error. On the left-hand side, we fix the number of
multipoles to be exploited in the reconstruction of the kernel; on the right-hand side, we consider a sort of "oracle" estimator, where the number of
multipoles grows with the optimal rate $N^{\frac{1}{2\beta _{\ast }-1}}$. \
As before, we took $C_{\ell ;Z}\simeq const\times \ell ^{-2}$, $\phi _{\ell
}\simeq const \times \ell ^{-\beta }$ for $\beta =2,2.5,3$; for $N=100,300,700$ the left -hand side uses $L_{N}\sim N^{0.6}$, while the right-hand side takes $L_{N}\sim N^{\frac{1}{%
2\beta _{\ast }-1}}$, as explained above.

We note how the estimators perform very efficiently, and show the errors
scale approximately as $N^{\alpha }$, where ${\alpha }\approx \frac{2\beta
_{\ast }-2}{2\beta _{\ast }-1}$, as predicted by our computations, see
Remark \ref{optimal_choice}. In particular, for\ $\beta _{\ast }=2$ our
results predict an upper bound for the $L^{2}$ error decaying as $N^{-0.67},$
whereas simulations show a decay in the order of $N^{-0.66};$ for $\beta
_{\ast }=2.5$ we have $N^{-0.75}$ and $N^{-0.82},$ and finally for $\beta
_{\ast }=3$ the predicted upper bound is in the order of $N^{-0.80}$, while the observed decay is in the order of $N^{-0.92}.$

\newpage

\begin{center}
\vspace{.5cm}
\begin{tabular}{rrrr} 
$N$ & Variance & Bias & MSE \\ \hline\hline
100 & 0.00174 & 0.00001 & 0.00175 \\
300 & 0.00116 & 0.00000 & 0.00117 \\
700 & 0.00084 & 0.00000 & 0.00084 \\ \hline
\end{tabular}%
\qquad
\begin{tabular}{rrrr}
$N$ & Variance & Bias & MSE \\ \hline\hline
100 & 0.00041 & 0.00023 & 0.00065 \\
300 & 0.00021 & 0.00010 & 0.00031 \\
700 & 0.00012 & 0.00005 & 0.00018 \\ \hline
\end{tabular}
\captionof{table}{$L^2$ errors obtained with $\beta_* = 2$; $L_N \sim N^{0.6}$ (left) and $L_N \sim N^{\frac{1}{2\beta_* - 1}}$ (right).}\label{Table1}

\vspace{.5cm}
\begin{tabular}{rrrr} 
$N$ & Variance & Bias & MSE \\ \hline\hline
100 & 0.00172 & 0.00001 & 0.00172 \\
300 & 0.00117 & 0.00000 & 0.00117 \\
700 & 0.00084 & 0.00000 & 0.00084 \\ \hline
\end{tabular}%
\qquad
\begin{tabular}{rrrr}
$N$ & Variance & Bias & MSE \\ \hline\hline
100 & 0.00028 & 0.00133 & 0.00162 \\
300 & 0.00013 & 0.00052 & 0.00065 \\
700 & 0.00008 & 0.00025 & 0.00033 \\ \hline
\end{tabular}
\captionof{table}{$L^2$ errors obtained with $\beta_* = 2.5$; $L_N \sim N^{0.6}$ (left) and $L_N \sim N^{\frac{1}{2\beta_* - 1}}$ (right).}\label{Table2}

\vspace{.5cm}
\begin{tabular}{rrrr}
$N$ & Variance & Bias & MSE \\ \hline\hline
100 & 0.00174 & 0.00000 & 0.00174 \\
300 & 0.00116 & 0.00000 & 0.00116 \\
700 & 0.00084 & 0.00000 & 0.00084 \\ \hline
\end{tabular}
\qquad
\begin{tabular}{rrrr}
$N$ & Variance & Bias & MSE \\ \hline\hline
100 & 0.00017 & 0.00440 & 0.00457 \\
300 & 0.00009 & 0.00072 & 0.00082 \\
700 & 0.00004 & 0.00072 & 0.00076 \\ \hline
\end{tabular}
\captionof{table}{$L^2$ errors obtained with $\beta_* = 3$; $L_N \sim N^{0.6}$ (left) and $L_N \sim N^{\frac{1}{2\beta_* - 1}}$ (right).}\label{Table3}
\vspace{.3cm}
\end{center}

We can now focus quickly on the main result of our paper, dealing with the Quantitative Central Limit Theorem, in Wasserstein distance; the latter is computed following the Python routine
(\texttt{scipy.stats.wasserstein\underline{ }distance}). We consider again a model where the autoregressive parameter and the angular power spectra are exactly the same as in the previous settings, in particular taking $\beta=3$ and $d=0.3,0.4$; we fix $L_{\max}=100$ for the number of components under the null hypothesis. Under these circumstances, we evaluate (univariate) Wasserstein distances for the kernel estimators at $m=9$ different locations, performing $B=1000$ Monte Carlo replications.

The results are reported in Table \ref{Table4}; here we take $B=1000$ Monte
Carlo replications, and taking $N=10^2,10^3,10^4$ for the cardinality of the time-domain observations. It should be noted that huge
sample sizes are quite common when dealing with sphere cross time data, see,
e.g., the NCEP/NCAR reanalysis datasets \cite{NCEP} for athmospheric
research. \\

Again, we note as simulations track closely the theoretical predictions. More precisely, for the theoretical upper bounds we expect $d_W$ to decay as $N^{\frac{1}{2}+d(1-\beta_{*})}$, leading to $N^{-0.1}$ in the setting of Table \ref{Table4}, $N^{-0.3}$ for Table \ref{Table5}, whereas estimates from the simulations give as worst rates $N^{-0.18}$ and $N^{-0.3}$, respectively.


\begin{center}
\vspace{.3cm}
\begin{tabular}{r|rrrrrrrrr}
$N \backslash z$ & -0.8 & -0.6 & -0.4 & -0.2 & 0 & 0.2 & 0.4 & 0.6 & 0.8 \\
\hline\hline
100 & 0.55 & 2.36 & 1.43 & 0.56 & 2.67 & 2.87 & 0.94 & 2.67 & 4.34 \\
1000 & 0.16 & 0.80 & 0.58 & 0.59 & 0.84 & 0.56 & 1.22 & 0.08 & 1.37 \\
10000 & 0.24 & 0.03 & 0.19 & 0.23 & 0.28 & 0.36 & 0.44 & 0.55 & 0.15 \\
\hline
\end{tabular}
\captionof{table}{Wasserstein distances obtained with $\beta_*=3$ and $L_N \sim N^{0.3}$.}\label{Table4}

\vspace{.5cm}
\begin{tabular}{r|rrrrrrrrr} 
$N \backslash z$ & -0.8 & -0.6 & -0.4 & -0.2 & 0 & 0.2 & 0.4 & 0.6 & 0.8 \\
\hline\hline
100 & 0.39 & 0.31 & 0.17 & 0.36 & 0.28 & 0.38 & 0.34 & 0.64 & 0.12 \\
1000 & 0.06 & 0.06 & 0.06 & 0.08 & 0.09 & 0.15 & 0.18 & 0.21 & 0.10 \\
10000 & 0.03 & 0.03 & 0.03 & 0.04 & 0.03 & 0.03 & 0.04 & 0.06 & 0.03 \\
\hline
\end{tabular}
\captionof{table}{Wasserstein distances obtained with $\beta_*=3$ and $L_N \sim N^{0.4}$.}\label{Table5}
\vspace{.3cm}
\end{center}

Clearly a full assessment of these procedures would require a much deeper numerical investigation; these preliminary results, however, seem rather encouraging for future developments.

\newpage

  \begin{minipage}{\linewidth}
      \centering
      \begin{minipage}{0.45\linewidth}
          \begin{figure}[H]
           \includegraphics[width=\linewidth]{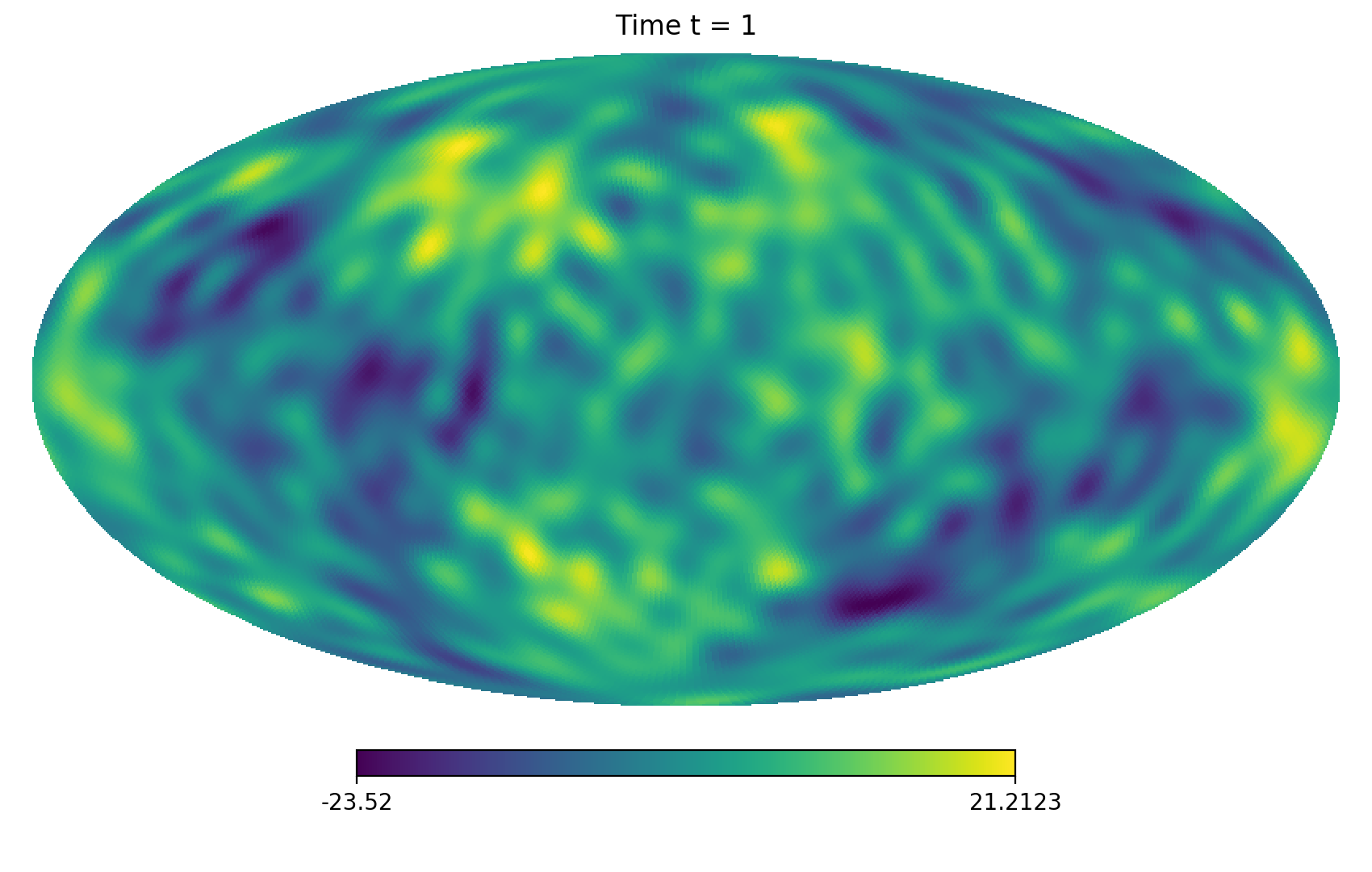}  \\
               \includegraphics[width=\linewidth]{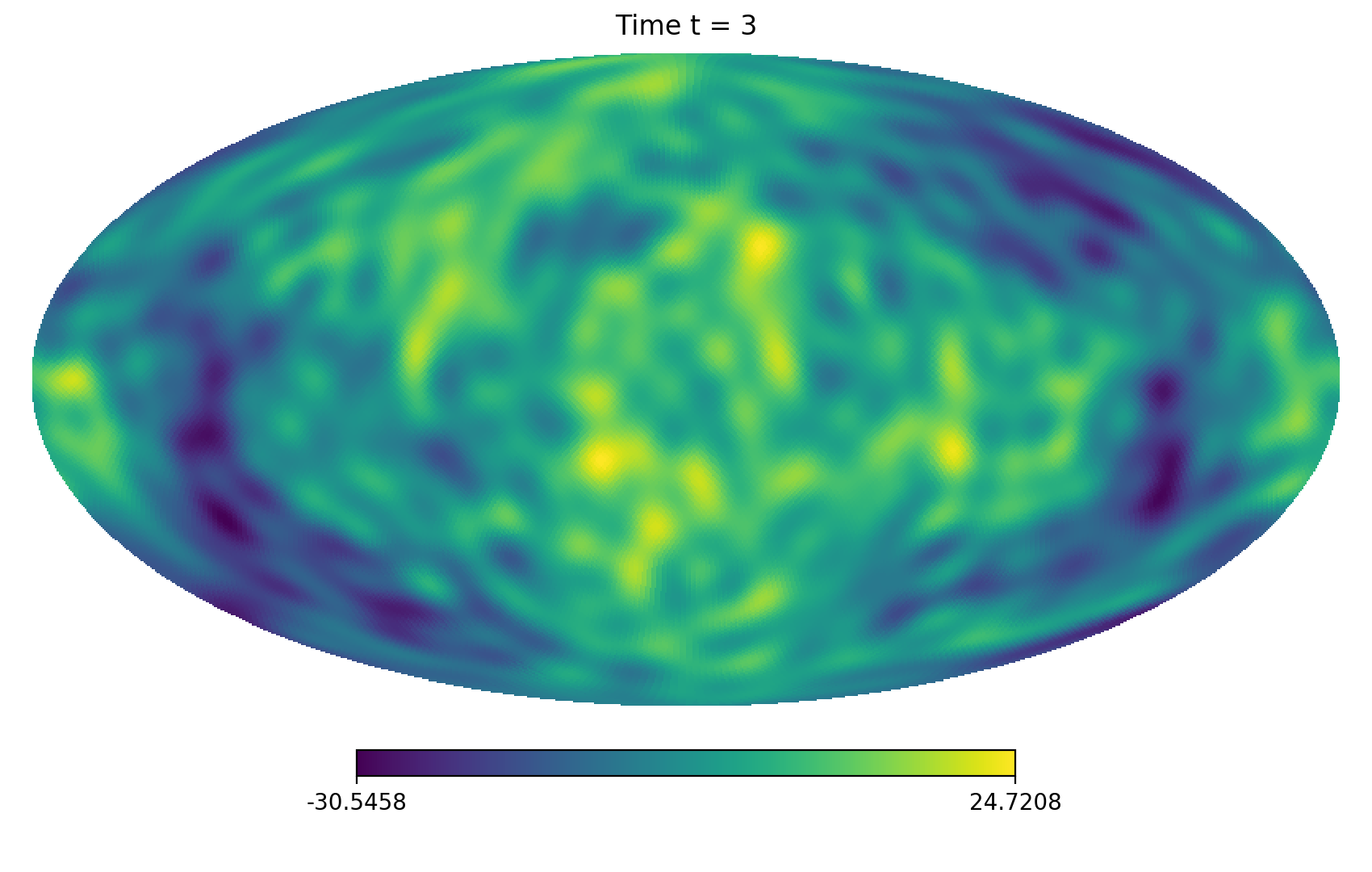}
          \end{figure}
      \end{minipage}
      \hspace{0.05\linewidth}
      \begin{minipage}{0.45\linewidth}
          \begin{figure}[H]
          	 \includegraphics[width=\linewidth]{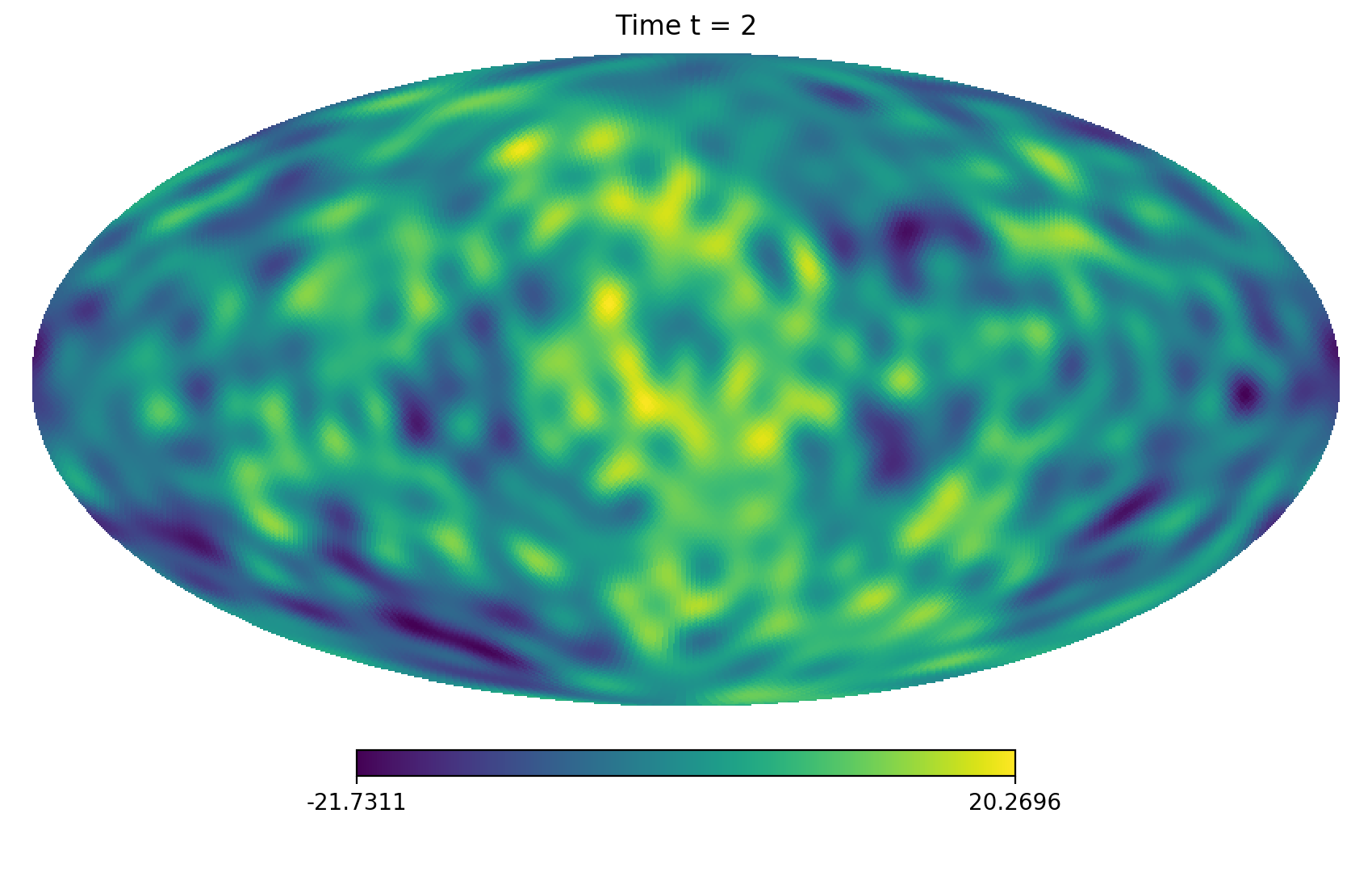}
               \includegraphics[width=\linewidth]{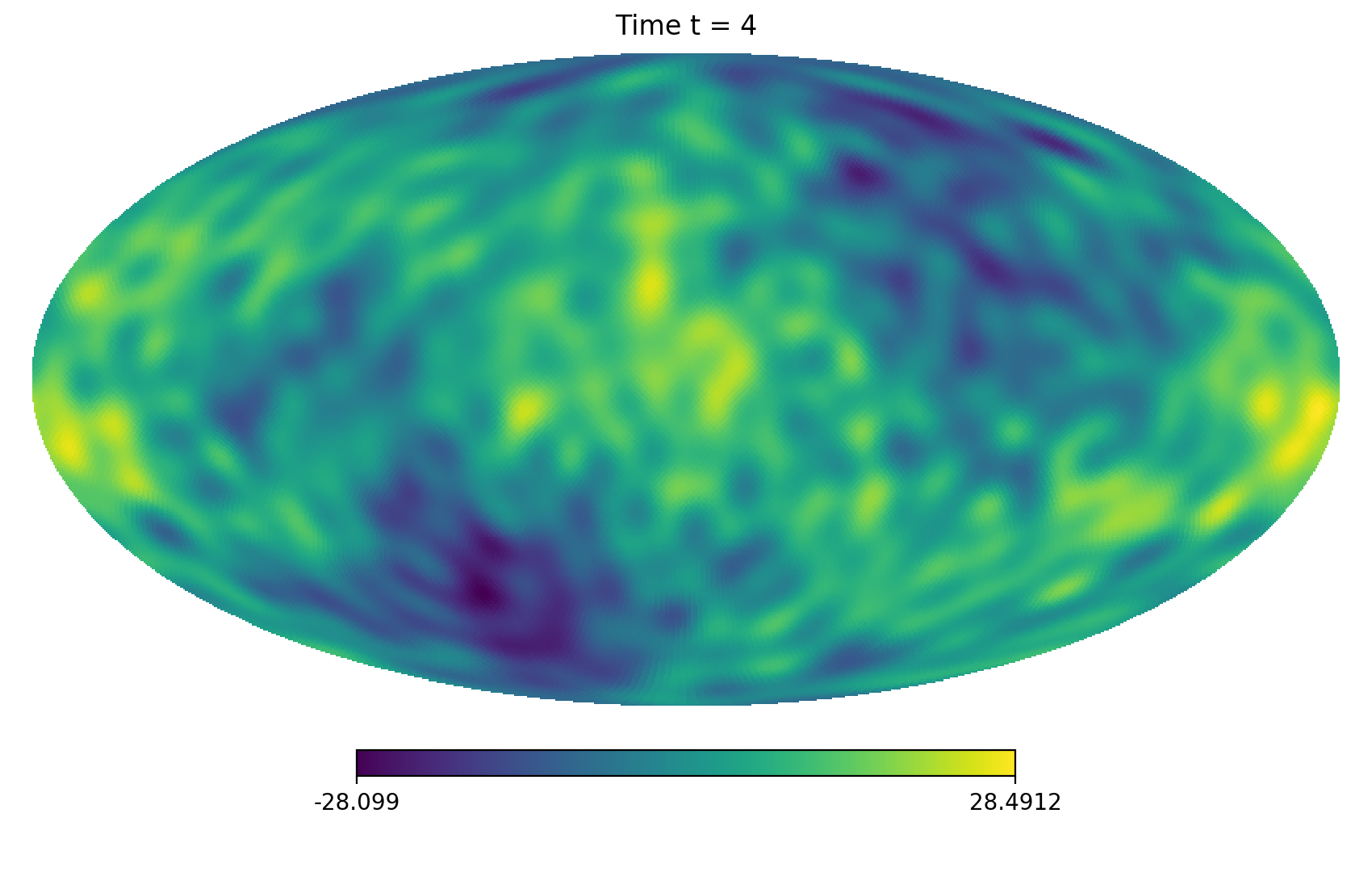}
          \end{figure}
      \end{minipage}
  \end{minipage}

  \begin{minipage}{\linewidth}
      \centering
      \begin{minipage}{0.45\linewidth}
          \begin{figure}[H]
           \includegraphics[width=\linewidth]{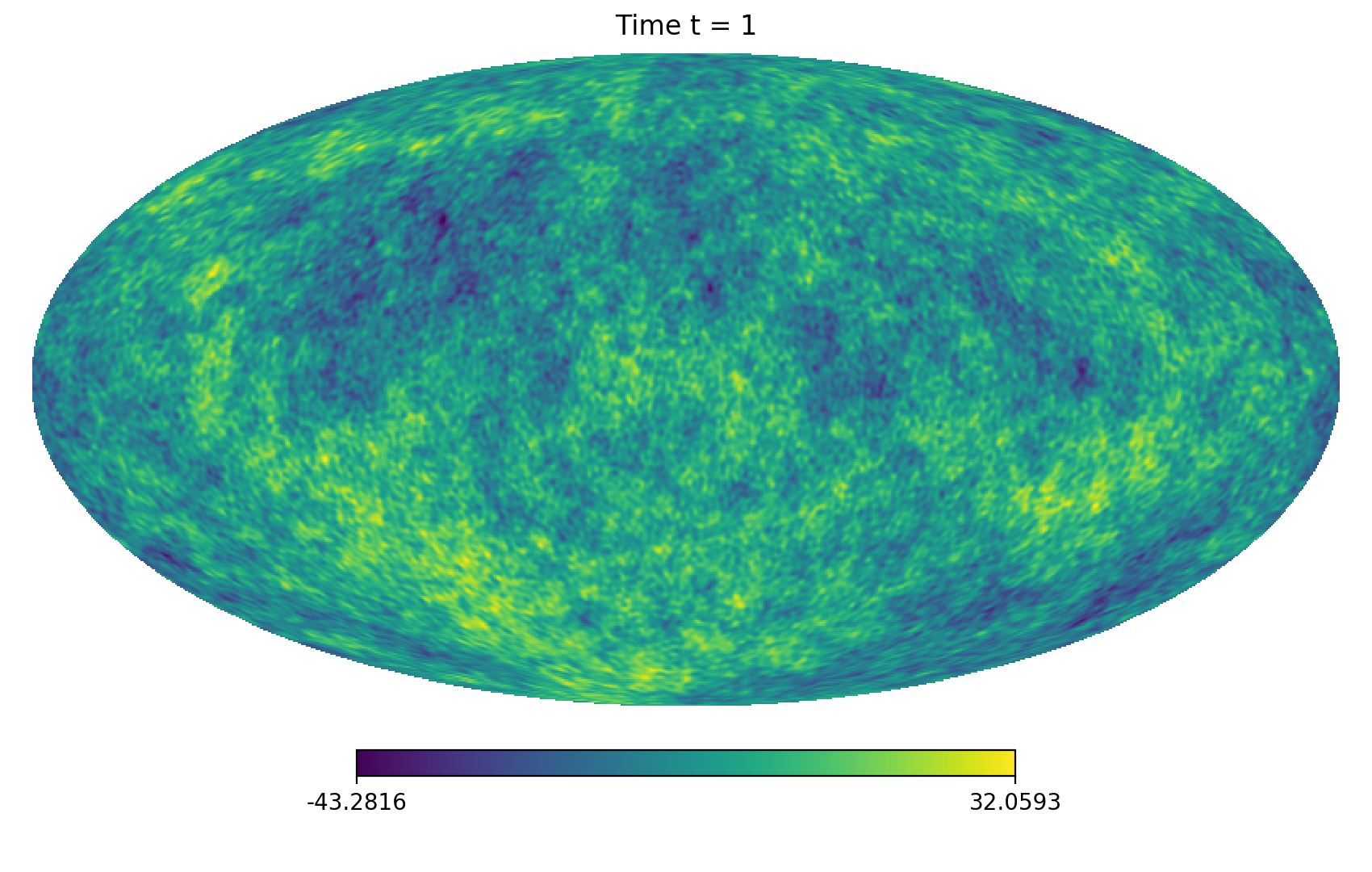}  \\
               \includegraphics[width=\linewidth]{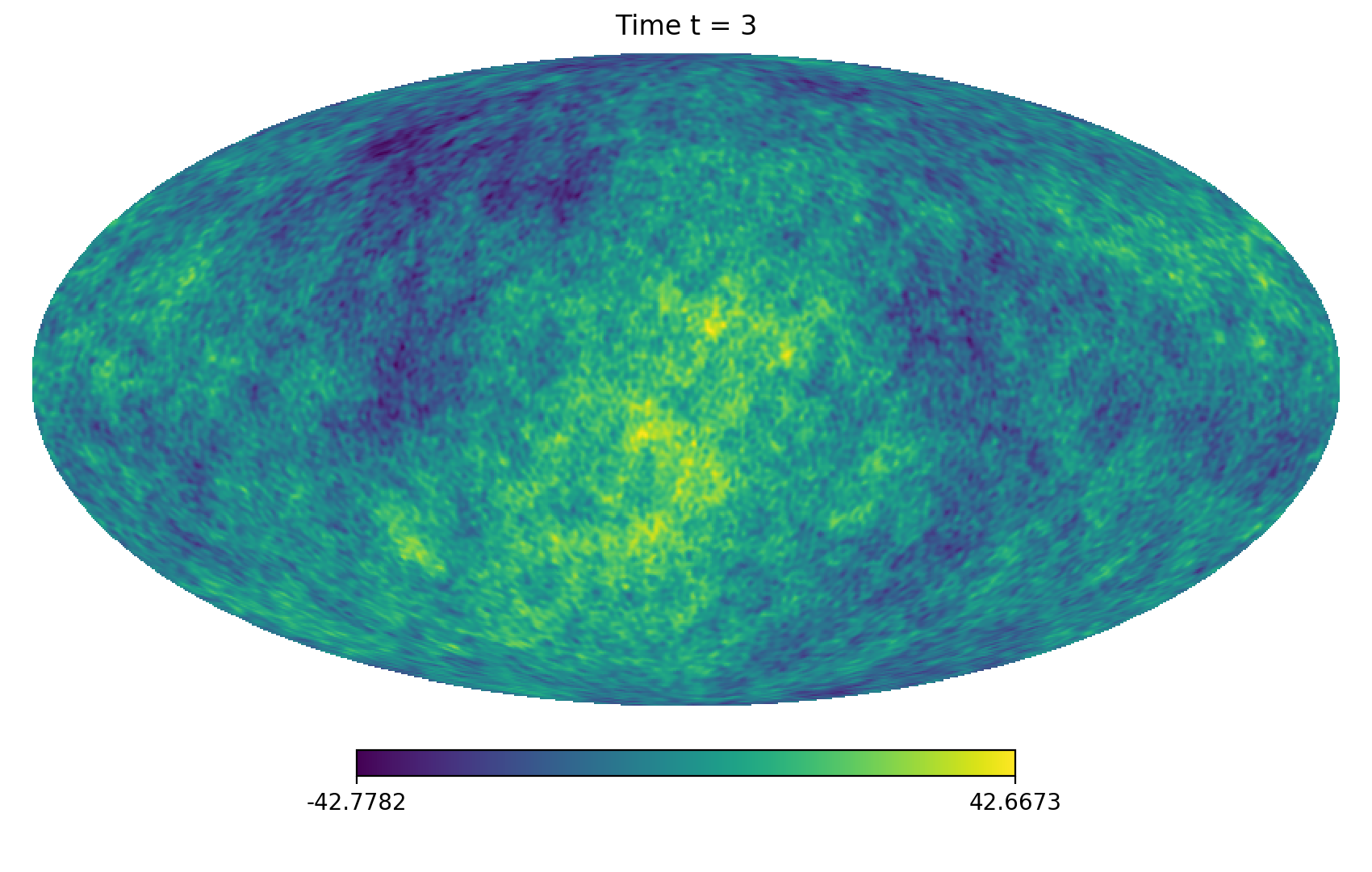}
          \end{figure}
      \end{minipage}
      \hspace{0.05\linewidth}
      \begin{minipage}{0.45\linewidth}
          \begin{figure}[H]
          	 \includegraphics[width=\linewidth]{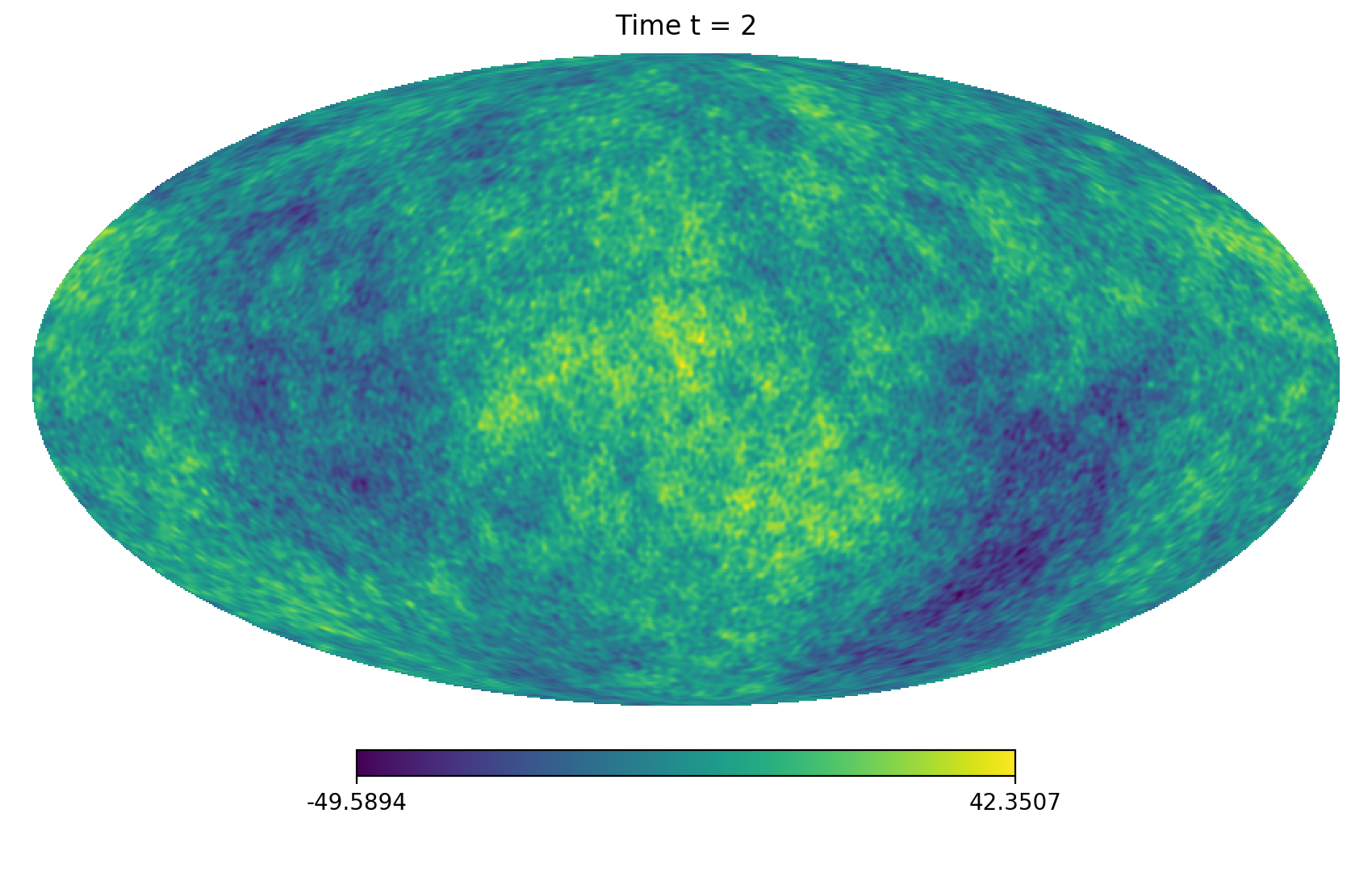}
               \includegraphics[width=\linewidth]{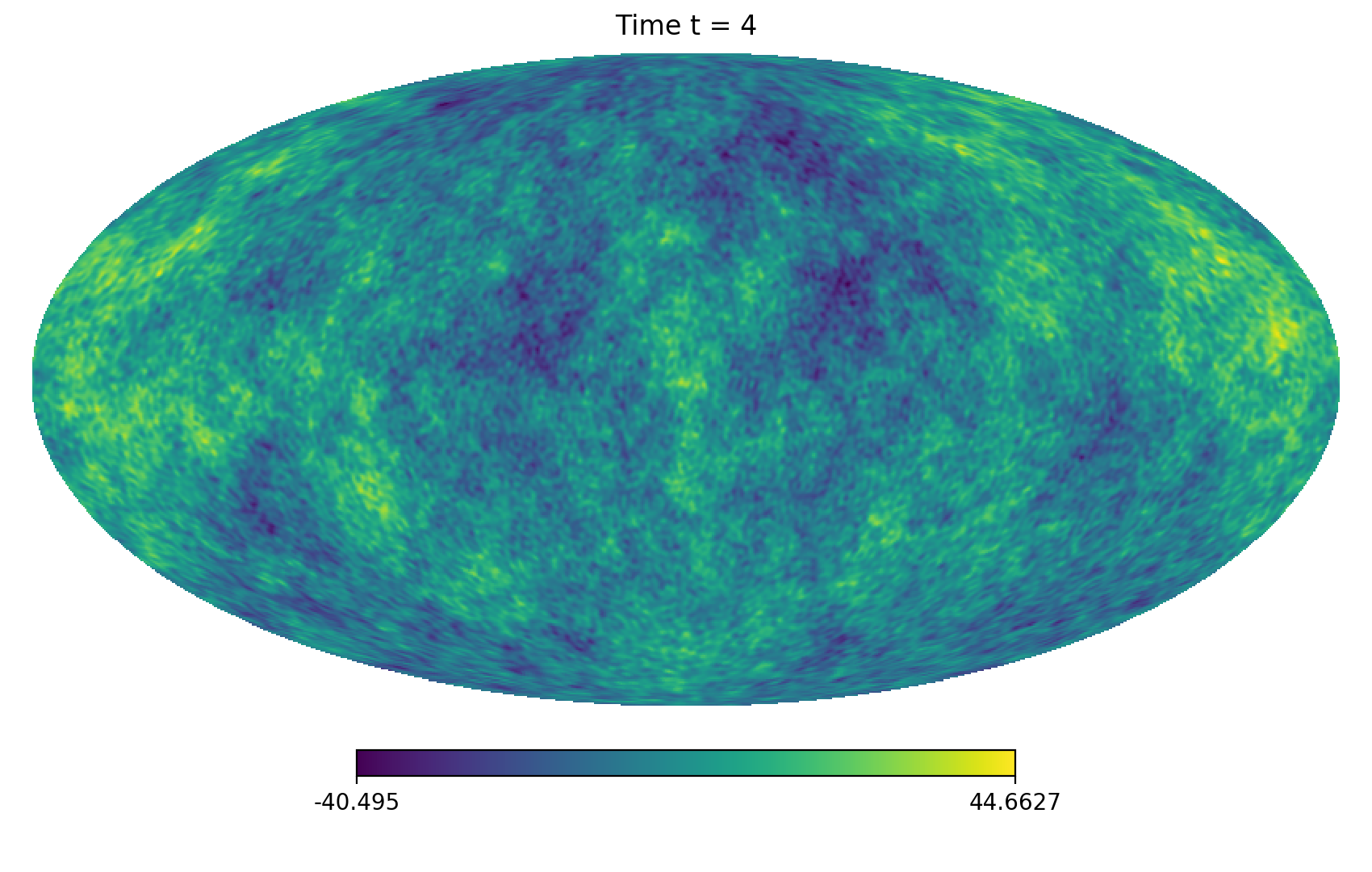}
          \end{figure}
      \end{minipage}
  \end{minipage}
        \captionof{figure}{Two realizations of sphere cross time random fields at time $t=1,2,3,4$. Upper panel: maximum resolution $L_{\max}=30$. Lower panel: maximum resolution $L_{\max}=200$.}

\newpage

\begin{frontmatter}

\title{Supplement to "Asymptotics for Spherical Functional Autoregressions"}
\runtitle{Supplement: Spherical Functional Autoregressions}

%
\author{\fnms{Alessia} \snm{Caponera}\ead[label=e1]{alessia.caponera@uniroma1.it}\thanksref{t1} and
\fnms{Domenico} \snm{Marinucci}\ead[label=e2]{marinucc@mat.uniroma2.it}\thanksref{t1}}

\address{Piazzale Aldo Moro, 5\\ 00185 Roma \\Italy\\ \href{mailto:alessia.caponera@uniroma1.it}{E-mail: \textnormal{alessia.caponera@uniroma1.it}} 
}
\affiliation{Department of Statistical Sciences, Sapienza University of Rome}


\thankstext{t1}{DM acknowledges the MIUR Excellence Department Project awarded to the Department of Mathematics, University of Rome Tor Vergata, CUP E83C18000100006. We are also grateful to Pierpaolo Brutti for many insightful suggestions and conversations.}
\address{Via della Ricerca Scientifica, 1\\00133 Roma \\Italy\\\href{mailto:marinucc@mat.uniroma2.it}{E-mail: \textnormal{marinucc@mat.uniroma2.it}}
}
\affiliation{Department of Mathematics, University of Rome Tor Vergata}
\runauthor{Caponera and Marinucci}




\end{frontmatter}


\setcounter{section}{0}
\renewcommand*{\theHsection}{chY.\the\value{section}}

\setcounter{equation}{0}
\renewcommand*{\theHequation}{chY.\the\value{equation}}

\setcounter{theorem}{0}
\renewcommand*{\theHtheorem}{chY.\the\value{theorem}}

\section{Appendix}  \label{Appendix}

Throught this Appendix, we assume that Conditions 8 and 13 hold. Under these assumptions
the proof that Equation (7)
admits a unique stationary and isotropic solution
can be given along the same lines as in \cite{Bosq} and it is omitted for
brevity's sake; see \cite{Caponera} for more discussion and details. Note
that, under these two Conditions,
the variance $C_{\ell }$ can be written in terms of the coefficients $\phi _{\ell ;j}
$, $j=1,\dots,p$, the autocorrelations $\rho _{\ell }(j)=C_{\ell
}(j)/C_{\ell }$, $j=1,\dots ,p$, and the error variance $C_{\ell;Z}$; namely
\begin{equation*}
C_{\ell }=\frac{C_{\ell ;Z}}{1-\phi _{\ell;1 }\rho _{\ell}(1)-\cdots -\phi
_{\ell; p}\rho _{\ell }(p)}>0\ , \qquad \ell \ge 0 \ .
\end{equation*}%
Hence,
\begin{equation*}
0<\frac{C_{\ell ;Z}}{C_{\ell }}=1-\phi _{\ell;1 }\rho _{\ell}(1)-\cdots
-\phi _{\ell; p}\rho _{\ell }(p) \ ,
\end{equation*}%
and there exists a positive constant $\phi ^{\ast }$ such that, uniformly over $\ell$,
\begin{equation}
\sum_{j=1}^{p}\phi _{\ell;j }\rho _{\ell }(j)\leq \phi ^{\ast }<1 \ . 
\label{phi_starstar_p}
\end{equation}%
Recall that $C_{\ell ;Z}/C_{\ell }$ and $C_{\ell }/C_{\ell ;Z}$ are (in
absolute value) bounded by positive constants since both converge to 1 as $%
\ell \rightarrow \infty $. Now, we denote with $g_{\ell }(\lambda )$ the correlation spectral
density
\begin{align*}
g_{\ell }(\lambda )& :=\frac{f_{\ell }(\lambda )}{C_{\ell }}=\frac{1}{2\pi }%
\sum_{\tau =-\infty }^{\infty }\rho _{\ell }(\tau )e^{i\lambda \tau } \\
& =\frac{1}{2\pi }\frac{1-\phi _{\ell;1 }\rho _{\ell}(1)-\cdots -\phi
_{\ell; p}\rho _{\ell }(p)}{|1-\phi _{\ell;1 }e^{i\lambda }-\cdots -\phi
_{\ell;p }e^{i\lambda p}|^{2}}\ ,\qquad \lambda \in \lbrack -\pi ,\pi ] \ ,
\end{align*}%
where $\rho_\ell(\cdot):=C_\ell(\cdot)/C_\ell$ is the autocorrelation function, and we recall that $\Sigma _{\ell }$ is the $p\times p$
matrix of autocorrelations, with $ij$-th element $\rho _{\ell }(i-j)$.
Since $g_{\ell}(\cdot)$ is a continuous symmetric function on $[-\pi ,\pi
]$, it follows that (see \cite{xiao2012covariance})
\begin{equation}  \label{lambda-minmax}
2\pi \underline{g_{\ell}}\leq \lambda _{\min }(\Sigma _{\ell })\leq \lambda
_{\max }(\Sigma _{\ell })\leq 2\pi \overline{g_{\ell}}\ ,
\end{equation}%
where $\underline{g_{\ell}}$ and $\overline{g_{\ell}}$ are the minimum and
maximum of $g_{\ell}(\cdot )$ in $[-\pi ,\pi ]$, respectively; $\lambda
_{\min }(\Sigma _{\ell })$ and $\lambda _{\max }(\Sigma _{\ell })$ are the
minimum and maximum eigenvalues of $\Sigma _{\ell }$, respectively.
Moreover, because we assumed $g_{\ell}(\lambda )>0,\ \forall \lambda \in
\lbrack -\pi ,\pi ]$, from \eqref{lambda-minmax} we conclude that the
minimum eigenvalue is strictly positive (and hence bounded away from zero)
and $\Sigma _{\ell }$ is positive definite (and then invertible). Since $%
\Sigma _{\ell }$ is a $p\times p$ real symmetric positive definite matrix,
then
\begin{equation*}
\Vert \Sigma _{\ell }\Vert _{\text{op}}\ =\lambda _{\max }(\Sigma _{\ell
})\leq \text{Tr}(\Sigma _{\ell })=p\ , \quad \Vert \Sigma _{\ell }^{-1}\Vert
_{\text{op}}=\frac{1}{\lambda _{\min }(\Sigma _{\ell })}\leq \frac{1}{2\pi
\underline{g_{\ell}}}\text{ },
\end{equation*}%
\begin{equation*}
\text{ Tr}(\Sigma _{\ell }^{-1})\leq p\Vert \Sigma _{\ell }^{-1}\Vert _{%
\text{op}}\leq \frac{p}{2\pi \underline{g_{\ell}}}\ ,
\end{equation*}
where $\Vert A\Vert _{\text{op}}\ =\sqrt{\lambda _{\max }(A^{\prime }A)}$,
and $\text{Tr}(A)$ is the trace of $A$. In addition,
\begin{align*}
\frac{1}{2\pi \underline{g_{\ell}}}& =\max_{\lambda \in \lbrack -\pi ,\pi ]}%
\frac{1}{2\pi g_{\ell}(\lambda )}=\max_{\lambda \in \lbrack -\pi ,\pi ]}%
\frac{|1-\phi _{\ell ; 1}e^{i\lambda }-\cdots -\phi _{\ell;p }e^{i\lambda
p}|^{2}}{1-\phi _{\ell;1 }\rho _{\ell }(1)-\cdots -\phi _{\ell ; p}\rho
_{\ell }(p)} \\
& \leq \frac{const}{1-\sum_{j=1}^{p}\phi _{\ell;j }\rho _{\ell }(j)}\ ,
\end{align*}%
since
\begin{equation*}
|1-\phi _{\ell;1 }e^{i\lambda }-\cdots -\phi _{\ell;p }e^{i\lambda p}|\leq
1+\sum_{j=1}^{p}|\phi _{\ell;j }| \le const \ .
\end{equation*}%
Then, from Equation \eqref{phi_starstar_p}, we can conclude that, uniformly
over $\ell$,
\begin{equation*}
\frac{1}{2\pi \underline{g_{\ell}}}\leq C\ , \quad  \text{some } C>0 \text{ .}
\end{equation*}

The first result below will be exploited to prove convergence in probability
of the denominator for our estimators, while the second part gives the
fourth-cumulant bound which is crucial for Stein-Malliavin arguments. We recall here for convenience the equalities
\[
\sqrt{N(2\ell+1)}\left( \widehat{\boldsymbol{\phi}}_{\ell;N}-\boldsymbol{\phi}_{\ell} \right)= A_{\ell;N}^{-1}B_{\ell;N}\ ,
\]
where
\[
 A_{\ell;N}=\frac{1}{C_{\ell}N(2\ell+1)}X_{\ell;N}'X_{\ell;N}\ , \qquad \Sigma_{\ell}=\mathbb{E}A_{\ell;N}\ ,
\]
and
\[
B_{\ell;N}=\frac{1}{C_{\ell}\sqrt{N(2\ell+1)}}X_{\ell;N}'\boldsymbol{\varepsilon}_{\ell;N} \text{ .}
\]
\begin{lemma}
\label{lemma1:arp}For any integers $\ell \geq 0$, $N>p$, there exists $M>0$
such that
\begin{equation}
{\mathbb{E}} \left[ a_{\ell ;N}(i,j)-{\mathbb{E}}[a_{\ell ;N}(i,j)]\right]
^{2}\leq \frac{M}{N(2\ell +1)}\ ,\qquad i,j\in \{1,\dots ,p\}\text{ ,}
\label{l2-limit:arp}
\end{equation}%
and
\begin{equation}  \label{cum4::eq}
\textnormal{Cum}_{4}\left [\widetilde{b}_{\ell ;N}(i)\right]=\frac{6}{%
N(2\ell +1)}\left( \frac{C_{\ell ;Z}}{C_{\ell }}s_{\ell }(i,i)\right) ^{2}%
\text{ ,}
\end{equation}%
where $\widetilde{b}_{\ell ;N}(i)=\sum s_{\ell }(i,j)b_{\ell ;N}(j)$ is the $%
i$-th element of the $p$-dimensional vector $\widetilde{B}_{\ell ;N}=\Sigma
_{\ell }^{-1}{B}_{\ell ;N},$ and $s_{\ell }(i,j)$ are the elements of the
inverse matrix $\Sigma _{\ell }^{-1}$\ .
\end{lemma}

The following result shows that replacing $A_{\ell ;N}$ by its expected value $\Sigma_{\ell}$ in the definition of the OLS-like estimator $\widehat{\phi }_{\ell ;N}$ does not have any asymptotic effect, as $N\rightarrow
\infty $.

\begin{lemma}
\label{lemma3:arp}For any integers $\ell \geq 0$ and $N>7+p$, there exists
generic positive constants such that
\begin{equation}  \label{mse:ols}
{\mathbb{E}}\left\Vert \widehat{\boldsymbol{\phi} }_{\ell ;N}-\boldsymbol{%
\phi} _{\ell } \right\Vert ^{2}\leq \frac{const}{N(2\ell +1)}\text{ ,}
\end{equation}%
and%
\begin{equation}
{\mathbb{E}}\left\Vert \sqrt{N(2\ell +1)}\left( \widehat{\boldsymbol{\phi} }%
_{\ell ;N}-\boldsymbol{\phi} _{\ell }\right) -\Sigma _{\ell }^{-1}{B}_{\ell
;N}\right\Vert \leq \frac{const}{\sqrt{N(2\ell +1)}}\text{ .}
\label{l1-limit:arp}
\end{equation}
\end{lemma}

The following results entail that $\lim_{N\rightarrow \infty }V_{N}=I_{mp}$;
actually the Propositions below give also a uniform rate of convergence.

\begin{lemma}
\label{lemma::sigma_rate} If $\|\boldsymbol{\phi}_\ell\| \le \frac{\gamma}{\ell^{\beta}}, \ \beta > 1, \ \ell > 0$,
\begin{equation*}
\left \| \frac{C_{\ell;Z}}{C_\ell}\Sigma_\ell^{-1} - I_p \right\|_{\infty} =
\mathcal{O}\left ( \frac{1}{\ell^{\beta}}\right) \ .
\end{equation*}
\end{lemma}

The next result is technical; given the huge amount of work which has taken
place on Legendre polynomials, we expect that the statement could be known
already, but we failed to locate a reference and therefore we report a full
proof for the sake of completeness.

\begin{lemma}
\label{hilbs2} Let $z=\cos \theta ,\ \theta \in (0,\pi )$,
\begin{equation}
\lim_{L\rightarrow \infty }\frac{1}{L+1}\sum_{\ell =0}^{L}(2\ell +1)P_{\ell
}^{2}(\cos \theta )=\frac{2}{\pi \sin \theta }\text{ ;}  \label{hilbs}
\end{equation}%
on the other hand, for $\theta ,\theta ^{\prime }\in (0,\pi ),\ \theta \neq
\theta ^{\prime },$ as $L\rightarrow \infty $,
\begin{equation}
\frac{1}{L+1}\sum_{\ell =0}^{L}(2\ell +1)P_{\ell }(\cos \theta )P_{\ell
}(\cos \theta ^{\prime })=O\left (\frac{\log L}{L} \right)\text{ }.
\label{hilbs2::eq}
\end{equation}
\end{lemma}

We can now start with the proof of these Lemmas.

\begin{proof}[Proof (Lemma \ref{lemma1:arp})]
Let us start by
observing that the $ij$-th element of $A_{\ell;N}$, denoted by $a_{\ell;N}(i,j)$, has expected value
\begin{align*}
\Ex \left[ a_{\ell;N}(i,j)  \right ] &=  \Ex \left[  \frac{1}{N (2\ell +1)C_\ell} \sum_{t=p+1}^{n}  \summ a_{\ell, m}(t-i)a_{\ell, m}(t-j) \right] \nonumber \\
&=  \frac{1}{N (2\ell +1)C_\ell} \sum_{t=p+1}^{n}  \summ \Ex \left[ a_{\ell, m}(t-i)a_{\ell, m}(t-j) \right]  \nonumber\\
&= \rho_\ell(i-j)\ .
\end{align*}
Now, we have
\[
{\mathbb{E}}   \left[ a_{\ell;N}(i,j) - {\mathbb{E}} [ a_{\ell;N}(i,j)] \right]^2=  \sum_{tt'} \sum_{m m'} \frac{\Cov[a_{\ell, m} (t_1-i) a_{\ell, m}(t_1-j), a_{\ell} (t_2-i)a_{\ell, m}(t_2-j)]}{N^2(2\ell+1)^2 C_\ell^2}
\]
\[\overset{\tau= t_1-t_2}{=}  \frac{1}{N(2\ell+1)} \sum_{\tau=1-N}^{N-1} \left ( 1 - \frac{|\tau|}{N} \right ) \left [ \left ( \frac{C_\ell(\tau)}{C_\ell} \right)^2 + \frac{C_\ell(\tau + i-j)}{C_\ell} \frac{C_\ell(\tau +j -i)}{C_\ell} \right]\ .
\]
Now observe $\rho_\ell^2(\cdot)=(C_\ell(\cdot)/C_\ell)^2$, the squared autocorrelation function of the process, is nonnegative and summable; that is, there exists $\rho_\ell^* \in \mathbb{R}^+$ so that
$
\sum_{\tau=-\infty}^{+\infty}\rho^2_\ell(\tau)= \rho_\ell^* < \infty,
$
and
\begin{equation*}
\sum_{\tau=-\infty}^{\infty} \left | \rho_\ell(\tau + i-j) \rho_\ell(\tau +j -i) \right | \le  \sum_{\tau=-\infty}^{+\infty}\rho_\ell^2(\tau)\ ,
\end{equation*}
in view of the Cauchy-Schwartz inequality. Thus, it holds that
\begin{align*}
&\sum_{\tau=1-N}^{N-1} \left ( 1 - \frac{|\tau|}{N} \right )\left [ \left ( \frac{C_\ell(\tau)}{C_\ell} \right)^2 + \frac{C_\ell(\tau + i-j)}{C_\ell} \frac{C_\ell(\tau +j -i)}{C_\ell} \right] \\
\le& \sum_{\tau=1-N}^{N-1} \rho_\ell^2(\tau) + \sum_{\tau=1-N}^{N-1} |\rho_\ell(\tau + i-j) \rho_\ell(\tau +j -i)|  \nonumber \\
\le&  \sum_{\tau=-\infty}^{+\infty}\rho_\ell^2(\tau) +  \sum_{\tau=-\infty}^{+\infty}\rho_\ell^2(\tau) \nonumber \\
=& \, 2 \rho_\ell^*\ .
\end{align*}
On the other hand,
\begin{equation*}
\rho_\ell^* = 2\pi \int_{-\pi}^{\pi} [g_{\ell}(\lambda)]^2 d \lambda,
\end{equation*}
and
\begin{align*}
g_{\ell}(\lambda) &= \frac{1}{2\pi} \frac{1-\phi_{\ell;1}\rho_\ell(1)-\cdots- \phi_{\ell;p} \rho_\ell(p)}{|1-\phi_{\ell;1}e^{i\lambda}-\cdots- \phi_{\ell;p} e^{i\lambda p}|^2} \nonumber \\
&\le \frac{1}{2\pi} \frac{const}{\left (1-\xi^{-1}_{\ast}  \right )^{2p}}\ ,
\end{align*}
since
\begin{equation*}
1- \sum_{j=1}^p \phi_{\ell;j} \rho_\ell(j) \le 1+ \sum_{j=1}^p |\phi_{\ell;j}| \le const\ ,
\end{equation*}
and
\begin{equation*}
|1-\phi_{\ell;1}e^{i\lambda}-\cdots-\phi_{\ell;p} e^{i\lambda p}|  = \prod_{j=1}^p |1-\xi^{-1}_{\ell; j} e^{i\lambda}| \ge \prod_{j=1}^p( 1-|\xi^{-1}_{\ell; j}| ) \ge ( 1-\xi^{-1}_{\ast} )^p >0\ ,
\end{equation*}
see also \cite{Caponera}.
Then, $\rho_\ell^* \le  const$, uniformly over $\ell$.

In conclusion, uniformly over $\ell$ and $N$,
\begin{equation*}
{\mathbb{E}}   \left[ a_{\ell;N}(i,j) - {\mathbb{E}} [ a_{\ell;N} (i,j)]  \right ]^2
 \le  \frac{M}{N(2\ell+1)}\ , \qquad i, j \in \{1,\dots,p\}\ ,
\end{equation*}
$M>0$.

Let us now focus on the elements of $\widetilde{B}_{\ell;N} = \Sigma_\ell^{-1}{B}_{\ell;N}$;
they are given by
\begin{equation*}
\widetilde{b}_{\ell;N}(i) = \sum_{j=1}^p s_\ell(i,j)b_{\ell;N}(j)\ , \qquad i =1,\dots, p \ .
\end{equation*}
These elements can be shown to satisfy the following properties:
\begin{enumerate}
\item $\Ex \left [\widetilde{b}_{\ell;N}(j)\right ] = \sum_{j=1}^p s_{\ell}(i,j) \Ex [b_{\ell;N}(j)]$ = 0; \label{B:7}
\item $\Ex \left [ \widetilde{b}_{\ell;N}(i) \widetilde{b}_{\ell;N}(j) \right]= s_{\ell}(i,j)\frac{C_{\ell; Z}}{C_\ell}$, \label{B:8}
since
\begin{equation*}
\Ex [ \Sigma_\ell^{-1}{B}_{\ell;N} (\Sigma_\ell^{-1}{B}_{\ell;N})' ]  = \Sigma_\ell^{-1} \Ex [ {B}_{\ell;N} {B}'_{\ell;N} ] \Sigma_\ell^{-1} = \frac{C_{\ell; Z}}{C_\ell} \Sigma_\ell^{-1}\ ,
\end{equation*}
and because
\begin{align*}
\Ex [ b_{\ell;N}(i) b_{\ell;N}(j)]&=\frac{1}{C_\ell^2} \frac{1}{N(2\ell+1)} \sum_{t t'} \sum_{mm'} \Ex[ a_{\ell, m}(t-i) a_{\ell, m; Z} (t)a_{\ell, m}(t'-j) a_{\ell, m; Z} (t')]  \nonumber\\
&= \frac{1}{C_\ell^2}  \frac{1}{N(2\ell+1)} \sum_{t t'} \sum_m \Ex \left [a_{\ell, m}(t-i) a_{\ell, m; Z} (t)a_{\ell, m}(t'-j) a_{\ell, m; Z} (t')\right]  \nonumber \\
&= \frac{1}{C_\ell^2}\frac{1}{N(2\ell+1)}  \sum_{t m} C_\ell(i-j)C_{\ell; Z} \nonumber\\
&= \frac{C_\ell(i-j)C_{\ell; Z}}{C_\ell^2}\ .
\end{align*}
\item $\text{Cum}_4\left [\widetilde{b}_{\ell;N}(i)\right ]= \frac{6}{N(2\ell+1)} \left (s_\ell(i,i) \frac{C_{\ell; Z}}{C_\ell} \right)^2$. \label{B:9}
\end{enumerate}

To compute $\text{Cum}_4\left [\widetilde{b}_{\ell;N}(i)\right ]$ we use once again the multilinearity property of cumulants, the real expansion and the diagram formula, so that we obtain:
\begin{equation*}
\text{Cum}_4[\widetilde{b}_{\ell;N}(i) ]= \sum_{j_1 j_2 j_3 j_4} s_\ell(i,j_1)s_\ell(i,j_2) s_\ell(i,j_3)s_\ell(i,j_4) \text{Cum}[b_{\ell;N}(j_1), b_{\ell;N}(j_2), b_{\ell;N}(j_3), b_{\ell;N}(j_4)],
\end{equation*}
with $\text{Cum}[b_{\ell;N}(j_1), b_{\ell;N}(j_2), b_{\ell;N}(j_3), b_{\ell;N}(j_4)]=\text{Cum}(j_1, j_2, j_3, j_4)$ given by
\begin{align*}
&\text{Cum}(j_1, j_2, j_3, j_4)= \frac{1}{C_\ell^4}\frac{1}{N^2(2\ell+1)^2}\\
\times&  \sum_{tm} \text{Cum} \big [ a_{\ell, m}(t-j_1) a_{\ell, m; Z} (t), a_{\ell, m}(t-j_2) a_{\ell, m; Z} (t),a_{\ell, m}(t-j_3) a_{\ell, m; Z} (t), a_{\ell, m}(t-j_4) a_{\ell, m; Z} (t)\big]   \nonumber\\
=&\, \frac{1}{N(2\ell+1)} \bigg[ 2 \frac{C_\ell(j_1 - j_2)}{C_\ell} \frac{C_\ell(j_3 - j_4)}{C_\ell} \left (\frac{C_{\ell; Z}}{C_\ell} \right)^2 \nonumber \\
+&\, 2 \frac{C_\ell(j_1 - j_3)}{C_\ell} \frac{C_\ell(j_2 - j_4)}{C_\ell} \left (\frac{C_{\ell; Z}}{C_\ell} \right)^2 \nonumber \\
+&\, 2 \frac{C_\ell(j_1 - j_4)}{C_\ell} \frac{C_\ell(j_2 - j_3)}{C_\ell} \left (\frac{C_{\ell; Z}}{C_\ell} \right)^2 \bigg]\ .
\end{align*}
Hence,
\begin{equation*}
\text{Cum}_4\left [\widetilde{b}_{\ell;N}(i) \right]= \frac{6}{N(2\ell+1)} \left (s_\ell(i,i) \frac{C_{\ell; Z}}{C_\ell} \right)^2,
\end{equation*}
as claimed.
\end{proof}

\begin{proof}[Proof (Lemma \ref{lemma3:arp})]
First, rewrite
\begin{equation*}
\sqrt{N(2\ell+1)}\left(\widehat{\boldsymbol{\phi}}_{\ell;N}  - \boldsymbol{\phi}_\ell \right)= \Sigma_\ell^{-1} {B}_{\ell;N} +  [ A_{\ell;N}^{-1} - \Sigma_{\ell}^{-1} ]{B}_{\ell;N}\ .
\end{equation*}

Since
\begin{align*}
\left \| \sqrt{N(2\ell+1)}\left(\widehat{\phi}_{\ell;N}  - \phi_\ell\right) -  \Sigma_\ell^{-1} {B}_{\ell;N} \right \| &= \|[ A_{\ell;N}^{-1} - \Sigma_{\ell}^{-1} ] {B}_{\ell;N} \| \nonumber \\
&=  \| [ I_p - \Sigma_{\ell}^{-1}A_{\ell;N}] A_{\ell;N}^{-1}  {B}_{\ell;N}  \| \nonumber \\
&\le  \|  I_p - \Sigma_{\ell}^{-1}A_{\ell;N}  \|_\op \|A_{\ell;N}^{-1} \|_\op \| {B}_{\ell;N} \|,
\end{align*}
we have
\begin{align*}
\Ex \| [ A_{\ell;N}^{-1} - \Sigma_{\ell}^{-1} ] {B}_{\ell;N} \|&\le \left ( \Ex \| A_{\ell;N}^{-1} \|^2_\op \| {B}_{\ell;N} \|^2 \right )^{1/2} \left ( \Ex \|I_p - \Sigma_{\ell}^{-1}A_{\ell;N} \|^2_\op \right)^{1/2} \nonumber\\
&\le \left ( \Ex \| A_{\ell;N}^{-1} \|^4_\op \right)^{1/4} \left(\Ex \| {B}_{\ell;N} \|^4 \right )^{1/4} \left ( \Ex \|I_p - \Sigma_{\ell}^{-1}A_{\ell;N} \|^2_\op \right)^{1/2},
\end{align*}
where
\begin{align*}
\Ex \| {B}_{\ell;N}  \|^4 &=  \sum_{i=1}^p \sum_{j=1}^p \Ex\left  [ b^2_{\ell;N}(i) b^2_{\ell;N}(j)\right ] \nonumber \\
&\le \sum_{i=1}^p \sum_{j=1}^p \left (\Ex\left  [ b^4_{\ell;N}(i) \right]  \right)^{1/2} \left (\Ex \left [ b^4_{\ell;N}(j)\right ] \right)^{1/2} \nonumber\\
&= \sum_{i=1}^p \sum_{j=1}^p  \left [  \frac{6}{N(2\ell +1)} \left ( \frac{C_{\ell; Z}}{C_\ell} \right)^2 + 3 \left ( \frac{C_{\ell; Z}}{C_\ell} \right)^2 \right] \nonumber \\
&< p^2  \left (\frac{24}{N(2\ell +1)} + 12 \right)\ ,
\end{align*}
and, from \eqref{l2-limit:arp},
\begin{align*}
\Ex \|I_p - \Sigma_{\ell}^{-1}A_{\ell;N} \|^2_\op &\le \| \Sigma_{\ell}^{-1}\|^2_\op\Ex \| \Sigma_\ell - A_{\ell;N}  \|^2_\op \nonumber \\
&\le const \sum_{i=1}^p \sum_{j=1}^p \Ex \left[ a_{\ell;N}(i,j) -  \Ex [a_{\ell;N}(i,j) ] \right]^2 \nonumber \\
&\le \frac{const}{N(2\ell+1)} \ .
\end{align*}

By definition,
\begin{align*}
\| A_{\ell;N} ^{-1} \|_\op = N(2\ell+1)\frac{C_\ell}{\lambda_{\min}(X'_{\ell;N}X_{\ell;N})}\ .
 \end{align*}
Since $X'_{\ell;N}X_{\ell;N}$ is a real symmetric $p\times p$ matrix,
\begin{equation*}
\lambda_{\min}(X_{\ell;N}'X_{\ell;N}) = \min_{\|\boldsymbol{\gamma}\|=1} \boldsymbol{\gamma}' X_{\ell;N}'X_{\ell;N} \boldsymbol{\gamma}\ .
\end{equation*}
$X_{\ell;N}'X_{\ell;N}$ can be seen as the sum of $2\ell +1$ independent matrix, i.e.
\begin{equation*}
X_{\ell;N}'X_{\ell;N} = \summ X'_{\ell, m;N}X_{\ell, m;N}\ ,
\end{equation*}
where $X_{\ell, m;N}$ is a $N \times p$ matrix, defined by (recalling that $n=N+p$)
\begin{equation*}
X_{\ell, m;N}=
\begin{pmatrix}
a_{\ell, m}(p)& a_{\ell, m}(p+1)& \cdots& a_{\ell, m}(n-1) \\
\vdots& \vdots& \vdots& \vdots \\
a_{\ell, m}(1)& a_{\ell, m}(2)& \cdots& a_{\ell, m}(n-p)
\end{pmatrix}\ .
\end{equation*}
Then,
\begin{align}\label{lambda_min_XX}
\lambda_{\min}(X'_{\ell;N}X_{\ell;N})&= \min_{\|\boldsymbol{\gamma}\|=1} \boldsymbol{\gamma}'  \left[\summ X'_{\ell, m;N}X_{\ell, m;N}\right] \boldsymbol{\gamma} \nonumber \\
&= \min_{\|\boldsymbol{\gamma}\|=1} \summ \boldsymbol{\gamma}' X'_{\ell, m;N}X_{\ell, m;N} \boldsymbol{\gamma}\ .
\end{align}

Now recall that $\Sigma_{\ell}$ is the $p\times p$ matrix of autocorrelations;
similarly we define $\Sigma_{\ell;N}$ as the $N\times N$ matrix of autocorrelations. Both are invertible since we assumed that the spectral density
\begin{equation*}
g_{\ell}(\lambda) = \frac{1}{2\pi} \frac{1-\phi_{\ell;1}\rho_\ell(1)-\cdots-\phi_{\ell;p} \rho_\ell(p)}{|1-\phi_{\ell;1}e^{i\lambda}-\cdots- \phi_{\ell;p}e^{i\lambda p}|^2}\ , \qquad \lambda \in [-\pi, \pi]\ ,
\end{equation*}
is a continuous positive function.

$X_{\ell, m;N}$ is a zero-mean Gaussian matrix with $\Ex[ X_{\ell, m;N} X'_{\ell, m;N}]= p C_\ell \Sigma_{\ell;N}$ and $\Ex[ X'_{\ell, m;N}X_{\ell, m;N} ]= N C_\ell \Sigma_{\ell}$, therefore it can be written as $X_{\ell, m;N} = (C_\ell \Sigma_{\ell;N})^{1/2} Z_{\ell, m;N}$, where $Z_{\ell, m;N}$ is a zero-mean Gaussian matrix with independent rows. If $\Sigma_{\ell;N}=P\Lambda P'$, where $P$ is an orthogonal matrix of eigenvectors and $\Lambda$ is the diagonal matrix of eigenvalues, then
\begin{align*}
\boldsymbol{\gamma}' X'_{\ell, m;N}X_{\ell, m;N}\boldsymbol{\gamma} &= C_\ell  \boldsymbol{\gamma}' Z'_{\ell, m;N} \Sigma_{\ell;N} Z_{\ell, m;N}\boldsymbol{\gamma} \nonumber \\
&= C_\ell \boldsymbol{\gamma}' Z'_{\ell, m;N} P\Lambda P' Z_{\ell, m;N}\boldsymbol{\gamma} \nonumber \\
&\ge C_\ell \lambda_{\min}(\Sigma_{\ell;N}) \boldsymbol{\gamma}' Z'_{\ell, m;N} P P' Z_{\ell, m;N} \boldsymbol{\gamma} \nonumber \\
&= C_\ell \lambda_{\min}(\Sigma_{\ell;N}) \boldsymbol{\gamma}' Z'_{\ell, m;N} Z_{\ell, m;N}\boldsymbol{\gamma}\ ,
\end{align*}
where $Z'_{\ell, m;N} Z_{\ell, m;N}$ is a Wishart random matrix with $N$ degrees of freedom. The same argument applies to all $2\ell + 1$ components of \eqref{lambda_min_XX}, so that
\begin{align}\label{lambda_min_XX_2}
\frac{\lambda_{\min}(X'_{\ell;N}X_{\ell;N})}{C_\ell}&\ge \lambda_{\min}(\Sigma_{\ell;N}) \min_{\|\boldsymbol{\gamma}\|=1} \summ \boldsymbol{\gamma}' Z'_{\ell, m;N} Z_{\ell, m;N}\boldsymbol{\gamma} \nonumber\\
&=\lambda_{\min}(\Sigma_{\ell;N}) \min_{\|\boldsymbol{\gamma}\|=1} \boldsymbol{\gamma}'  \left[ \summ Z'_{\ell, m;N} Z_{\ell, m;N}\right] \boldsymbol{\gamma}  \nonumber \\
&=  \lambda_{\min}(\Sigma_{\ell;N}) \min_{\|\boldsymbol{\gamma}\|=1} \boldsymbol{\gamma}' Z'_{\ell;N} Z_{\ell;N} \boldsymbol{\gamma}' \nonumber \\
&=  \lambda_{\min}(\Sigma_{\ell;N}) \lambda_{\min}(Z'_{\ell;N}Z_{\ell;N})\ .
\end{align}
The summation in \eqref{lambda_min_XX_2} includes $2\ell + 1$ independent Wishart random matrix each with $N$ degrees of freedom and $\Sigma_\ell$ as scale matrix, then $Z'_{\ell;N} Z_{\ell;N}$ is a Wishart random matrix with $N(2\ell +1)$ degrees of freedom and $ \Sigma_\ell$ as scale matrix, and  $\lambda_{\min}(Z'_{\ell;N}Z_{\ell;N})$ its minimum eigenvalue. Furthermore, this result guarantees the invertibility of the matrix $X'_{\ell;N} X_{\ell;N}$.

By the standard inequality on trace and operator norms for matrices, we obtain that
\begin{align*}
\Ex\| A_{\ell;N} ^{-1} \|^4_\op &\le \frac{N^4(2\ell+1)^4}{\lambda^4_{\min}(\Sigma_{\ell;N})} \Ex\| (Z'_{\ell;N}Z_{\ell;N}) ^{-1} \|^4_\op\nonumber \\
&\le \frac{N^4(2\ell+1)^4}{\lambda^4_{\min}(\Sigma_{\ell;N})}   \Ex \left[  \tr(  (Z'_{\ell;N}Z_{\ell;N})^{-1})\right]^4 \nonumber\\
&\le \frac{N^4(2\ell+1)^4}{(2\pi \underline{g_{\ell}})^4}\Ex \left[  \tr(  (Z'_{\ell;N}Z_{\ell;N})^{-1})\right]^4\ .
\end{align*}
For $N(2\ell +1) > 7 + p$ the fourth moment of the trace of an inverse Wishart matrix is given in \cite{Matsumoto}:
\begin{align*}
u_4(\eta)\Ex \left[  \tr(  (Z'_{\ell;N}Z_{\ell;N})^{-1})\right]^4 &= 48(5\eta -3) \tr(\Sigma_\ell^{-4}) \nonumber\\
&+128\eta(\eta-2) \tr(\Sigma_\ell^{-3}) \tr(\Sigma_\ell^{-1}) \nonumber \\
&+12(2\eta^2-5\eta+9)(\tr(\Sigma_\ell^{-2}))^2 \nonumber \\
&+ 12(4\eta^3-12\eta^2+3\eta+3)\tr(\Sigma_\ell^{-2}) (\tr(\Sigma_\ell^{-1}))^2 \nonumber\\
&+ (\eta+1) (2\eta-3)(4\eta^2-12\eta+1)(\tr(\Sigma_\ell^{-1}))^4\ ,
\end{align*}
where $\eta = \frac{N(2\ell+1)}{2} - \frac{p+1}{2}$, and
\begin{equation*}
u_4(\eta)= 2^4\eta(\eta-1)(\eta-2)(\eta-3)(2\eta-1)(\eta+1)(2\eta+1)(2\eta+3)\ .
\end{equation*}
If $\lambda_{\ell;1},\dots,\lambda_{\ell;p}$ are the eigenvalues of $\Sigma_\ell$, we have
\begin{equation*}
0 < \tr(\Sigma_\ell^{-k}) = \sum_{j=1}^p\left ( \frac{1}{\lambda_{\ell;j}}\right)^k \le \left(  \sum_{j=1}^p \frac{1}{\lambda_{\ell;j}}\right)^k = (\tr(\Sigma_\ell^{-1}))^k\ ,
\end{equation*}
$k \ge 1$. Then, for $2\eta > 7 + p$,
\begin{align*}
u_4(\eta)\Ex \left[  \tr(  (Z'_{\ell;N}Z_{\ell;N})^{-1})\right]^4 &\le (\tr(\Sigma_\ell^{-1}))^4(8\eta^4+20\eta^3+10\eta^2-5\eta-3)\nonumber \\
&=(\tr(\Sigma_\ell^{-1}))^4 (2\eta-1)(\eta+1)(2\eta+1)(2\eta+3)\ ,
\end{align*}
and
\begin{align*}
\Ex \left[  \tr(  (Z'_{\ell;N}Z_{\ell;N})^{-1})\right]^4  &\le \frac{(\tr(\Sigma_\ell^{-1}))^4}{2^4 \eta (\eta -1)(\eta-2)(\eta-3)}\nonumber \\
&= \frac{(\tr(\Sigma_\ell^{-1}))^4}{\prod_{k=1}^4 (N(2\ell+1) - p +1  - 2n ))}\ .
\end{align*}
In addition,
\begin{align*}
\frac{N^4(2\ell+ 1)^4}{\prod_{k=1}^4 \left(N(2\ell+1) - p +1 - 2n\right)} &= \frac{1}{\prod_{k=1}^4 \left(1 - \frac{p - 1 + 2n}{N(2\ell+1)}\right)} \nonumber \\
& \le \frac{1}{ \left(1 - \frac{p + 7}{N(2\ell+1)}\right)^4} \nonumber \\
&\le \frac{1}{ \left(1 - \frac{p + 7}{p + 8}\right)^4}\ ,
\end{align*}
for every $\ell \ge 0$ and $N > 7 + p$. Thus, \eqref{l1-limit:arp} holds.

The second part of this Lemma follows easily, indeed
\begin{align*}
\Ex\left \|\widehat{ \boldsymbol{\phi}}_{\ell;N} - \boldsymbol{\phi}_\ell \right \|^2 &= \frac{1}{N(2\ell+1)} \Ex\left \|\sqrt{N(2\ell+1)} ( \widehat{ \boldsymbol{\phi}}_{\ell;N} - \boldsymbol{\phi}_\ell )\right \|^2  \\
& =  \frac{1}{N(2\ell+1)} \Ex\left \| A^{-1}_{\ell;N} {B}_{\ell;N}\right \|^2  \\
& \le \frac{1}{N(2\ell+1)} \left (\Ex\left \| A^{-1}_{\ell;N}\right \| ^4 \right)^{1/2}\left( \Ex\left \|{B}_{\ell;N}\right \|^4 \right)^{1/2} \\
& \le \frac{const}{N(2\ell+1)}\ ,
\end{align*}
in view of the bounds that we just established on the fourth-moments of the norms of $A_{\ell,N}^{-1}$ and $B_{\ell,N}$.
\end{proof}

\begin{proof}[Proof (Lemma \ref{lemma::sigma_rate})]
We first need to prove that $\lim_{\ell \to \infty} \Sigma_\ell =  I_p$, where we recall that $\Sigma_\ell$ is the matrix of autocorrelations $\rho_\ell(i-j)$. For $i=j$, $\rho_\ell(i-j)=1$, for all $\ell$; on the other hand, for $i \ne j$,
\begin{equation*}
\rho_\ell(i-j)= \phi_{\ell;1} \rho_\ell( i-j - 1) + \cdots + \phi_{\ell;p}\rho_\ell( i-j - p)\ ,
\end{equation*}
and
\begin{equation*}
|\rho_\ell(i-j)| \le \sum_{k=1}^p |\phi_{\ell;k}| \to 0\ , \qquad \ell \to \infty\ .
\end{equation*}

For $\ell >0$,
\begin{equation*}
\left \| \frac{C_{\ell;Z}}{C_\ell}\Sigma_\ell^{-1} - I_p \right\|_{\infty}  \le \left \| I_p -\frac{C_{\ell}}{C_{\ell;Z}} \Sigma_\ell \right \|_{\infty} \left \|\frac{C_{\ell;Z}}{C_\ell} \Sigma_\ell^{-1} \right \|_{\infty} \le const \left | \frac{C_{\ell;Z}}{C_\ell} \right | \left \| I_p -\frac{C_{\ell}}{C_{\ell;Z}} \Sigma_\ell \right \|_{\infty}\ .
\end{equation*}
Moreover, since
\begin{equation*}
|\rho_\ell(i-j)| \le p\|\boldsymbol{\phi}_\ell\| \le \frac{p\gamma}{\ell^{\beta}} \ , \qquad i\ne j\ ,
\end{equation*}
and
\begin{equation*}
\left | 1 - \frac{C_\ell}{C_{\ell;Z}} \right | = \left|\frac{C_\ell}{C_{\ell;Z}}  \right |\left |\sum_{j=1}^p \phi_{\ell;j} \rho_\ell (j) \right| \le \left|\frac{C_\ell}{C_{\ell;Z}}  \right | p \| \boldsymbol{\phi}_\ell \| \le \left|\frac{C_\ell}{C_{\ell;Z}}  \right |  \frac{p\gamma}{\ell^{\beta}}\ ,
\end{equation*}
we have
\begin{equation*}
\left \| \frac{C_{\ell;Z}}{C_\ell}\Sigma_\ell^{-1} - I_p \right\|_{\infty}  \le \frac{const}{\ell^{\beta}}\ ,
\end{equation*}
as claimed.
\end{proof}

The last proof is for the technical Lemma on summation of squared Legendre polynomials.

\begin{proof}[Proof (Lemma \ref{hilbs2})]
For $\ell \ge 1$, by Hilb's asymptotics (see \cite{Szego},\cite{Wig}), it holds that
\begin{equation*}
P_\ell(\cos \theta) = \sqrt{\frac{2}{\pi \ell \sin \theta}} \sin \left ( \ell \theta + \alpha \right) + \mathcal{O}\left( \ell^{-3/2}\right)\ , \qquad 0 < \theta < \pi\ ,
\end{equation*}
with $\alpha= \frac{\theta}{2} +  \frac{\pi}{4}$.
Then,
\begin{align*}
(2\ell+1)P_\ell^2(\cos \theta) &=(2\ell+1)\left( \sqrt{\frac{2}{\pi \ell \sin \theta}} \sin \left ( \ell \theta + \alpha \right) + \mathcal{O}\left( \ell^{-3/2}\right) \right )^2 \nonumber \\
&= \frac{4}{\pi \sin \theta}  \sin^2 \left ( \ell \theta+ \alpha \right) +  \mathcal{O}\left( \ell^{-1}\right)\ , \qquad 0 < \theta < \pi\ .
\end{align*}
In view of the standard identities
\begin{equation*}
\sin x = \frac{e^{ix} - e^{-ix}}{2i}\ ,
\end{equation*}
and
\begin{equation*}
\sum_{k=0}^{n-1} e^{i x k} = \frac{1- e^{i x n}}{1-e^{ix}}\ , \qquad x \ne 0\ ,
\end{equation*}
we have
\begin{align*}
\sum_{\ell=1}^L \sin^2(\ell \theta + \alpha) &= \sum_{\ell=1}^L \left (  \frac{e^{i(\ell \theta + \alpha)} - e^{-i(\ell \theta + \alpha)}}{2i} \right )^2 \nonumber \\
&= - \frac{1}{4} \sum_{\ell=1}^L \left [ e^{i2(\ell \theta + \alpha)} + e^{-i2(\ell \theta + \alpha)} - 2 \right] \nonumber \\
&=- \frac{e^{i2(\theta+\alpha)}}{4} \left ( \frac{1- e^{i 2\theta L}}{1-e^{i2\theta}}\right) - \frac{e^{-i2(\theta+\alpha)}}{4} \left ( \frac{1- e^{-i 2\theta L}}{1-e^{-i2\theta}}\right) + \frac{1}{2}(L+1)\ ,
\end{align*}
hence,
\begin{equation*}
\lim_{L \to \infty} \frac{1}{L+1} \sum_{\ell=1}^L \sin^2(\ell \theta + \alpha)  = \frac{1}{2}\ .
\end{equation*}
Also, it holds that if $\lim_{k \to \infty} a_k = A$, $|A|< \infty$, then $\lim_{k \to \infty} \frac{1}{n} \sum_{k=1}^n a_k = A$. As a consequence,
\begin{equation*}
\lim_{L \to \infty } \frac{1}{L+1}\sum_{\ell=0}^L(2\ell+1)P_\ell^2(\cos \theta) = \frac{2}{ \pi \sin \theta}\ , \qquad \theta  \in (0,\pi)\ .
\end{equation*}

Likewise, for $\theta, \theta' \in (0, \pi), \ \theta \ne \theta'$,
\begin{align*}
(2\ell+1)P_\ell(\cos \theta) P_\ell(\cos \theta')&= \frac{4}{\pi \sqrt{\sin \theta \sin \theta'}}  \sin \left ( \ell \theta+ \alpha \right)  \sin \left ( \ell \theta'+ \alpha' \right) +  \mathcal{O}\left( \ell^{-1}\right)\ .
\end{align*}
As before, we have
\begin{align*}
\sum_{\ell=1}^L \sin(\ell \theta + \alpha) \sin(\ell \theta' + \alpha') &= \sum_{\ell=1}^L \left (  \frac{e^{i(\ell \theta + \alpha)} - e^{-i(\ell \theta + \alpha)}}{2i} \right ) \left (  \frac{e^{i(\ell \theta' + \alpha')} - e^{-i(\ell \theta' + \alpha')}}{2i} \right ) \nonumber \\
&= - \frac{1}{4} \sum_{\ell=1}^L \left [ e^{i(\ell (\theta +\theta') + \alpha + \alpha')} + e^{-i(\ell (\theta+\theta') + \alpha + \alpha')} \right] \nonumber \\
&+\frac{1}{4} \sum_{\ell=1}^L \left [ e^{i(\ell (\theta -\theta') + \alpha - \alpha')} + e^{-i(\ell (\theta-\theta') + \alpha - \alpha')}  \right] \nonumber \\
&=- \frac{e^{i(\theta+\theta'+\alpha+ \alpha')}}{4} \left ( \frac{1- e^{i (\theta+\theta')L}}{1-e^{i(\theta+\theta')}}\right) - \frac{e^{-i(\theta+\theta' + \alpha + \alpha')}}{4} \left ( \frac{1- e^{-i (\theta+\theta') L}}{1-e^{-i(\theta+\theta')}}\right) \nonumber \\
&+ \frac{e^{i(\theta-\theta'+\alpha - \alpha')}}{4} \left ( \frac{1- e^{i (\theta-\theta')L}}{1-e^{i(\theta-\theta')}}\right) + \frac{e^{-i(\theta-\theta'+\alpha - \alpha')}}{4} \left ( \frac{1- e^{-i (\theta-\theta') L}}{1-e^{-i(\theta-\theta')}}\right)\ ,
\end{align*}
hence,
\begin{equation*}
\frac{1}{L+1} \sum_{\ell=1}^L \sin(\ell \theta + \alpha) \sin(\ell \theta' + \alpha')=\mathcal{O} \left (\frac{1}{L} \right )\ .
\end{equation*}
In addition, since
$\sum_{\ell=1}^{L} \ell ^{-1} = \mathcal{O}(\log L)$,
we can then conclude that
\begin{equation*}
\frac{1}{L+1}\sum_{\ell=0}^L(2\ell+1)P_\ell(\cos \theta) P_\ell(\cos \theta')=\mathcal{O} \left (\frac{\log L}{L} \right)\ , \qquad \theta,\theta'  \in (0,\pi)\ , \ \theta \ne \theta'\ .
\end{equation*}
as $L \to \infty$.
\end{proof}

\begin{remark}
Note that (\ref{hilbs2::eq}) does not converge pointwise if $\theta $ or $%
\theta ^{\prime }=0;$ for instance, for $\theta =\theta ^{\prime }=0$ we
have $\frac{1}{L+1}\sum_{\ell =0}^{L}(2\ell +1)=L+1,$ whereas for $\theta
\neq 0$, $\theta ^{\prime }=0$ \ (\ref{hilbs2}) oscillates among given
constants.
\end{remark}

\end{document}